\documentclass[12pt]{amsart}

\usepackage[margin=1in]{geometry}
\usepackage{amsmath,amsthm,amssymb}
\usepackage{dsfont}
\usepackage{amsmath,amscd}
\usepackage{mathrsfs}
\usepackage{color}
\usepackage{hyperref}
\usepackage{enumitem}
\usepackage{array}
\usepackage{graphicx}
\usepackage{xcolor}
\usepackage{ulem}

\hypersetup{
    colorlinks=true, 
    linkcolor=blue, 
    urlcolor=red, 
    linktoc=all 
}

\input xy
\xyoption{all}

\newcommand{\N}{\mathbb{N}}

\newcommand{\C}{\mathbb{C}}

\newcommand{\PP}{\mathbb{P}}
\newcommand{\Z}{\mathbb{Z}}

\newcommand{\one}{\mathds{1}}

\newcommand{\im}{\textrm{Im}}

\DeclareMathOperator\coker{coker}

\theoremstyle{plain}
\newtheorem{thm}{Theorem}[section]
\newtheorem{prop}[thm]{Proposition}
\newtheorem{lemma}[thm]{Lemma}
\newtheorem{remark}[thm]{Remark}

\newtheorem{cor}[thm]{Corollary}

\theoremstyle{definition}

\newtheorem{eg}{Example}[thm]

\newcolumntype{C}[1]{>{\centering\let\newline\\\arraybackslash\hspace{0pt}}m{#1}}
\newcolumntype{L}[1]{>{\let\newline\\\arraybackslash\hspace{0pt}}m{#1}}

\begin{document}

\title[Suspension splittings of 5-dimensional Poincar\'{e} duality complexes]{Suspension splittings of 5-dimensional Poincar\'{e} duality complexes and their applications}
\author[S.~Amelotte]{Steven Amelotte}
\address{Department of Mathematics, University of Western Ontario, London ON, N6A 5B7, Canada} 
\email{samelot@uwo.ca}
\author[T.~Cutler]{Tyrone Cutler}
\address{Beijing Institute of Mathematical Sciences and Applications, Beijing 101408, China} 
\email{tyronecutler@bimsa.cn}
\author[T.~So]{Tseleung So}
\address{Department of Mathematics, University of Western Ontario, London ON, N6A 5B7, Canada} 
\email{tso28@uwo.ca}

\maketitle

\begin{abstract}
Let $X$ be a connected, orientable, 5-dimensional Poincar\'{e} duality complex with torsion-free $H_1(X;\Z)$. We show that $\Sigma X$ is homotopy equivalent to a wedge of recognisable spaces and study to what extent its homotopy type is determined by algebraic data. These results are then used to compute the unstable cohomotopy groups $\pi^3(X)$ and $\pi^3(X;\mathbb{Z}/k)$ as well as give partial information about the cohomotopy set $\pi^2(X)$.
\end{abstract}


\section{Introduction}

Smooth, simply connected, closed, orientable 5-manifolds are well-understood. In fact, the classification of these objects was started by Smale \cite{Smale:1962} and completed by Barden \cite{Barden:1965} in the 1960s (see also Zhubr \cite{Zhubr:2001} for another perspective). Later, St{\"o}cker \cite{Stocker:1982} would extend this to a classification of simply connected, 5-dimensional Poincar\'e duality complexes up to oriented homotopy type.

Comparatively little is known about non-simply connected 5-manifolds. Since every finitely-presentable group appears as the fundamental group of a smooth, closed, orientable~5-manifold \cite{ACampoKotschick:1994}, such manifolds are unclassifiable. Even the homotopy types of these manifolds are not well-understood.

However, results for non-simply connected 5-manifolds are starting to appear in the literature. These include the recent papers by Hambleton and Su \cite{HambletonSu:2013}, where certain 5-manifolds~$M$ with $\pi_1(M)=\mathbb{Z}/2$ are classified, and by Kreck and Su \cite{KreckSu:2017}, where certain~5-manifolds with free fundamental groups are considered. Despite this, there seems to be very little in the literature relating to non-simply connected Poincar\'e duality complexes in dimension 5.

In this paper we study the suspension splitting problem for connected, orientable, 5-dimensional Poincar\'{e} duality complexes $X$ for which $H_1(X;\Z)$ is torsion-free. The homology groups of such a space are given by
\begin{equation}\label{table_original M hmlgy}
\begin{tabular}{C{1.4cm}|C{1.4cm}|C{1.4cm}|C{1.4cm}|C{1.4cm}|C{1.4cm}|C{1.4cm}}
$i$	&$0$	&$1$	&$2$	&$3$	&$4$	&$5$\\
\hline
$H_i(X)$	&$\Z$	&$\Z^m$	&$\Z^n\oplus T$	&$\Z^{n}$	&$\Z^m$	&$\Z$
\end{tabular}
\end{equation}
where $T=\bigoplus^{\ell}_{i=1}\Z/t_i\Z$ is a torsion group. We assume each $t_i=p_i^{r_i}$ for some prime $p_i$. Let
\[
T=\{\{p_1^{r_1},\ldots,p_{\ell}^{r_{\ell}}\}\}
\]
be the collection containing all $t_i$'s, where repeated $t_i$'s are allowed.

In order to state the main theorem we need to set some notation. Firstly, for $n\geq2$ and~$k\geq1$ we let $P^n(k)=S^{n-1}\cup_ke^n$ be the Moore space obtained by attaching an $n$-cell to~$S^{n-1}$ by a degree $k$ map. Then it is known that (see~\cite{baues} or Lemma~\ref{lemma_hmtpy gps Moore space})
\[
\pi_4(P^3(2^r))\cong\begin{cases}
\Z/4\langle\xi_1\rangle	&\text{if }r=1\\
\Z/2\langle\xi_r\rangle\oplus\Z/2\langle\imath_r\circ\eta^2\rangle &\text{if }r>1
\end{cases}
\]
where $\eta^2\colon S^4\to S^2$ is the composite $S^4\overset{\eta}{\to}S^3\overset{\eta}{\to}S^2$ of Hopf maps, $\imath_r\colon S^2\to P^3(2^r)$ is the inclusion of the bottom cell, and $\xi_r$ is a lift of $\eta\colon S^4\to S^3$ through the pinch map $P^3(2^r)\to S^3$.
Secondly, let $\C\PP^2(2^r)$ be the mapping cone of the composite $S^3\overset{\eta}{\to}S^2\overset{\imath_r}{\to} P^3(2^r)$, and define~$\epsilon_r$ to be the composite
\[
\epsilon_r\colon S^4\overset{\xi_r}{\longrightarrow}P^3(2^r)\hookrightarrow\C\PP^2(2^r).
\]
In Lemma~\ref{lemma_pi_5 Sigma CP^2(2^r)} we show that the suspension of $\epsilon_r$ generates a $\mathbb{Z}/2$ summand in $\pi_5(\Sigma\C\PP^2(2^r))$.

Our main theorem states that $\Sigma X$ is homotopy equivalent to a wedge of recognisable spaces.
Denote the $d$-fold wedge sum of a space $A$ by $A^{\vee d}$.
Given a sequence $\{f_i\colon S^4\to A_i\}^d_{i=1}$ of maps, let $f_1\bot\cdots\bot f_d$ denote the composite
\begin{equation}\label{eqn_bot def}
f_1\bot\cdots\bot f_d\colon S^4\overset{\text{comult}}{\longrightarrow}\bigvee^d_{i=1}S^4\overset{\bigvee_if_i\,}{\longrightarrow}\bigvee^d_{i=1}A_i.
\end{equation}

\begin{thm}\label{thm_main theorem}
Let $X$ be a 5-dimensional orientable Poincar\'e duality complex with $H_1(X;\Z)$ torsion-free. Then there is a homotopy equivalence
\[
\Sigma X\simeq(S^2)^{\vee m}\vee(S^3)^{\vee n_3}\vee(S^4)^{\vee n_4}\vee(S^5)^{\vee n_5}\vee\bigvee_{t_i\in T'}P^4(t_i)\vee(\Sigma\C\PP^2)^{\vee b}\vee\bigvee^{c}_{j=1}\Sigma\C\PP^2(2^{r_j})\vee\Sigma C_f
\]
for some non-negative integers $n_3,n_4,n_5,b,c$, where $T'$ and $\{2^{r_j}\}$ are subcollections of $T$, and~$C_f$ is the mapping cone of one of the following maps:
\begin{enumerate}[label=\normalfont(\arabic*)]
\item\label{type of f_S^3 v Moore}
$x\eta\bot y{\xi}_{r}\bot z{\epsilon}_s\colon S^4\longrightarrow (S^3)^{\vee x}\vee P^3(2^r)^{\vee y}\vee\C\PP^2(2^s)^{\vee z}$
\item\label{type of f_S^2 v Moore}
$x\eta^2\bot y{\xi}_{r}\bot z{\epsilon}_s\colon S^4\longrightarrow (S^2)^{\vee x}\vee P^3(2^r)^{\vee y}\vee\C\PP^2(2^s)^{\vee z}$
\item\label{type of f_two Moore}
$x({\imath\circ\eta^2})\bot y{\xi}_{r}\bot z{\epsilon}_s \colon S^4\longrightarrow P^3(2^{q})^{\vee x}\vee P^3(2^r)^{\vee y}\vee\C\PP^2(2^s)^{\vee z}$
\item\label{type of f_one Moore}
$y({\imath\circ\eta^2}+{\xi}_{r})\bot z{\epsilon}_s\colon S^4\longrightarrow P^3(2^r)^{\vee y}\vee\C\PP^2(2^s)^{\vee z}$
\end{enumerate}
for some $x,y,z\in\{0,1\}$. Furthermore, if $y=z=1$ then $s<r$.
\end{thm}


These results complete a sequence of recent studies into suspension splittings of low-dimensional manifolds and Poincaré duality complexes. In particular, suspension splittings of 4-manifolds were considered by the third author and Theriault \cite{ST19}, and later by Li \cite{Li:2023}. Suspension splittings of 6-manifolds were studied by Huang~\cite{Huang:2023}, and also by the last two named authors~\cite{CS}. Suspension splittings of certain 7-manifolds are obtained in \cite{HuangLi:2023}, and for highly-connected manifolds of higher dimension in \cite{Huang:2022}. Other recent examples include results on the homotopy types of suspensions of toric manifolds~\cite{ChoiKajiTheriault:2017,HasuKishimotoSato:2016} and flag manifolds~\cite{KajiTheriault:2019}.

Suspension splittings of 5-manifolds were previously considered by Huang~\cite{Huang:2021}. However, these were not the main results of that paper, and Theorem~\ref{thm_main theorem} is sharper than the splittings found in~\cite{Huang:2021}.

We explain in Section~\ref{sect5} to what extent the homotopy type of $\Sigma X$ is determined by algebraic information contained in $H^*(X)$. An interesting case is when $X$ is a closed, orientable 5-manifold and we have access to its characteristic classes. All manifolds in this paper will be smooth, with Stiefel--Whitney classes $w_i(X)$ and Pontryagin class $p_1(X)$.
\begin{thm}\label{corollary1}
Let $X$ be a closed, orientable 5-manifold with torsion-free $H_1(X;\Z)$ and let~$C_f$ be the mapping cone described in Theorem~\ref{thm_main theorem}. Then the following statements hold.
\begin{enumerate}[label=\normalfont(\arabic*)]
\item $X$ is spin if and only if
\[
\Sigma X\simeq(S^2)^{\vee m}\vee(S^3)^{\vee n}\vee(S^4)^{\vee n}\vee(S^5)^{\vee m}\vee\bigvee_{t_i\in T}P^4(t_i)\vee S^6.
\]
\item If $X$ is nonspin and $w_3(X)\neq0$, then $f$ contains exactly one of $\xi_1$ or $\epsilon_1$. If $f$ contains~$\xi_1$, then it contains no $\epsilon_s$ for any  $s\geq1$. Furthermore, if $w_2(X)\cdot w_3(X)\neq0$, then $f$ contains $\xi_1$ if and only if $p_1(X)$ is divisible by $2$.
\item If $X$ is nonspin and $w_3(X)=0$, then $f$ must contain one of $\eta,{\xi}_r,{\epsilon}_s$. If either of~${\xi}_r$ or~${\epsilon}_s$ appears, then $r,s>1$.
The class $w_2(X)$ survives to exactly the $E_r$-page of the mod $2$ Bockstein spectral sequence if and only if either $\epsilon_r$ appears, or $\xi_r$ appears and no $\epsilon_s$ does for $s<r$. The class $w_2(X)$ survives to the $E_\infty$-page of the mod $2$ Bockstein spectral sequence if and only if $f\simeq\eta$.
\end{enumerate}
\end{thm}
Actually, every Poincaré duality complex has algebraically defined Stiefel-Whitney classes, and Theorem~\ref{corollary1} is a special case of more general results found in Propositions~\ref{propspin5complex},~\ref{w3nonzeroimpl}, and~\ref{w3zeroimplication}.

In addition to the theoretical interest of suspension splittings, they have many useful applications. In the present paper we develop two. Firstly, the splitting of Theorem~\ref{thm_main theorem} induces decompositions of $h^*(X)$ for any generalised cohomology theory $h^*(\;)$. Our result in this regard is stated in Section~\ref{sectgencohom} and is similar in spirit to the cohomological decomposition results for 4-manifolds given in \cite{ST19}, and for 6-manifolds given in \cite{CS}.

Our second set of applications is concerned with the calculation of {\it unstable} cohomotopy groups of $X$. In Section~\ref{sectcohomgroups} we describe how induced decompositions of the stable cohomotopy groups of $X$ can be leveraged to obtain information in the unstable range.
Cohomotopy groups of manifolds have been considered in many places in the literature. Along with Pontrjagin's classical work relating to 3-manifold cohomotopy groups, there are the more recent results of Kirby, Melvin, and Teichner~\cite{KirbyMelvieichner:2012}, and of Taylor~\cite{Taylor:2012} computing $\pi^2(X)$ for a 4-manifold $X$. More recently, the group $\pi^4(X)$ has been considered by Konstantis~\cite{Konstantis:2021} when~$X$ is a closed spin 5-manifold.

For a 5-manifold $X$, the most interesting cohomotopy group to study is $\pi^3(X)=[X,S^3]$, where the group structure is that induced by the Lie group structure of $S^3$. We also consider the cohomotopy sets with coefficients $\pi^q(X;\mathbb{Z}/k)=[X,P^{q+1}(k)]$. If $q\geq4$, then $\pi^q(X;\mathbb{Z}/k)$ has a canonical group structure, but otherwise $\pi^q(X;\mathbb{Z}/k)$ is a priori only a pointed set. In Propositions~\ref{proppi3X} and \ref{pi3withcoeffsisagroup} we obtain the following result.

\begin{thm}\label{labeltheorem2}
Let $X$ be a 5-dimensional CW complex. Then stabilisation induces a group isomorphism
\[
\pi^3(X)\cong \pi^3_S(X).
\]
Moreover, for any $k\geq1$, stabilisation induces a bijection
\[
\pi^3(X;\mathbb{Z}/k)\cong \pi^3_S(X;\mathbb{Z}/k)
\]
which equips $\pi^3(X;\mathbb{Z}/k)$ with a group structure. \qed
\end{thm}

Since stable cohomotopy (with $\mathbb{Z}$ or $\mathbb{Z}/k$ coefficients) is a generalised cohomology theory, both $\pi^*_S(X)$ and $\pi^\ast_S(X;\mathbb{Z}/k)$ are subject to the splitting result given in Theorem~\ref{Thm51}. Thus Theorem~\ref{labeltheorem2} leads to the following.

\begin{cor}\label{corolaboutpi3s}
If $X$ is a 5-dimensional orientable Poincar\'e duality complex with $H_1(X;\Z)$ torsion-free, then there are group isomorphisms
\begin{align*}
\pi^3(X)\cong& \, {\textstyle \pi_S^3(S^3)^{\oplus n_4}\oplus\pi_S^3(S^4)^{\oplus{n_5}}\oplus \bigoplus_{t'_i\in T'} \pi^3_S(P^3(t'_i))\oplus\pi_S^3(\mathbb{CP}^2)^{\oplus b}\oplus}\\
&\,{\textstyle \bigoplus^c_{j=1}\pi^3_S(\mathbb{CP}^2(2^{r_j}))\oplus \pi_S^3(C_f)}
\end{align*}
and
\begin{align*}
\pi^3(X;\mathbb{Z}/k)\cong&\;\textstyle{\pi^3_S(S^3;\mathbb{Z}/k)^{\oplus n_4}\oplus\pi^3_S(S^4;\mathbb{Z}/k)^{\oplus n_5}\oplus\bigoplus_{t'_i\in T'} \pi^3_S(P^3(t'_i);\mathbb{Z}/k)\oplus}\\
&\;\textstyle{\pi^3_S(\mathbb{CP}^2;\mathbb{Z}/k)^{\oplus b}\oplus\bigoplus^c_{j=1}\pi^3_S(\mathbb{CP}^2(2^{r_j});\mathbb{Z}/k)\oplus \pi^3_S(C_f;\mathbb{Z}/k).} \qed
\end{align*}
\end{cor}
Information on the summands appearing on the right-hand side of these isomorphisms is given in Sections~\ref{sectpi3X} and \ref{sectpi3xzk}. 

Finally, we give partial information about the set $\pi^2(X)$. Our results owe a heavy debt to Taylor's work \cite{Taylor:2012}, but appear to be new. In light of Theorem~\ref{labeltheorem2}, the next result implies that $\pi^2(X)$ is often determined by stable data.
\begin{thm}\label{taylorespi2xthrm}
Let $X$ be a closed, orientable, 5-dimensional manifold. Then the following statements hold.
\begin{enumerate}[label=\normalfont(\arabic*)]
\item If $X$ is simply connected, then there is a bijection
\[
\pi^2(X)\cong H^2(X)\times \pi^3(X).
\]
\item If $H_1(X;\Z)$ is torsion-free and $H_2(X;\Z)$ is torsion, then $\eta\colon S^3\rightarrow S^2$ induces a bijection
\[
\pi^2(X)\cong \pi^3(X).
\]
\end{enumerate}
\end{thm}
Further results regarding $\pi^2(X)$ are found in Section~\ref{sectpi2X}. We remark that a class of non-simply connected 5-manifolds satisfying the assumptions of Theorem~\ref{taylorespi2xthrm} (2) is studied in~\cite[Theorem 1.3]{KreckSu:2017}.

The paper is structured as follows. Section~\ref{sect2} contains a selection of lemmas as well as computations of certain homotopy groups which will be necessary in the sequel. In the short Section~\ref{sect3} we outline the proof of Theorem~\ref{thm_main theorem} and establish notation which will be used in the longer, more technical Section~\ref{section_homotopy types of W_4 and W_5}. This section is broken into two subsections which contain the details necessary to complete the proof of Theorem~\ref{thm_main theorem}. In Section~\ref{sect5} we analyse the homotopy type of $\Sigma X$ using Poincaré duality and prove Theorem~\ref{corollary1}.

Applications of Theorem~\ref{thm_main theorem} are given in Section~\ref{sect6}. Two subsections contain $(i)$ evaluation of generalised cohomology groups of~$X$, and $(ii)$ calculation of certain cohomotopy groups of~$X$. A proof of Theorem~\ref{taylorespi2xthrm} is found in this final subsection, as well as the details required to prove Theorem~\ref{labeltheorem2} and Corollary~\ref{corolaboutpi3s}.

\section*{Acknowledgements}
The second author would like to thank Stefan Behrens for useful discussion. The third author is supported by NSERC Discovery Grant and NSERC RGPIN-2020-06428.

\section{Preliminaries}\label{sect2}
Throughout the paper we assume that all spaces are path-connected and equipped with basepoints which are preserved by all maps. All homology and cohomology groups are reduced and taken with integer coefficients unless otherwise specified.

\subsection{Homotopy types of mapping cones}\label{sect_hmtpy type of mapping cone}

Given maps $f_i\colon S^n\to A_i$ for $1\leq i\leq a$, write~$C(f_1,\ldots,f_a)$ for the mapping cone of the map
\[
\sum^a_{i=1}(\jmath_i\circ f_i): S^n\to\bigvee^a_{i=1}A_i
\]
where $\jmath_k: A_k\hookrightarrow\bigvee^a_{i=1}A_i$ is the inclusion. Observe that in the notation introduced in equation ~\eqref{eqn_bot def} we have $\sum^a_{i=1}(\jmath_i\circ f_i)\simeq f_1\bot\cdots\bot f_a$. Thus $C(f_1,\ldots,f_a)$ is also the mapping cone of $f_1\bot\cdots\bot f_a$.

\begin{lemma}\label{lemma_simplify mapping cone}
Suppose each $A_i$ is the suspension of a connected CW complex and each map~$f_i$ is a suspension. Then for any map $g\colon A_k\to A_l$, there is a homotopy equivalence
\[
C(f_1,\ldots,f_a)\simeq C(f_1,\ldots,f_{l-1}, f_l+(g\circ f_k),f_{l+1},\ldots,f_a).
\]
\end{lemma}

\begin{proof}
Without loss of generality, we may assume that $k=1$ and $l=a=2$. Let $c: S^n\to S^n\vee S^n$ be the comultiplication. Since $f_{1}$ is a suspension, the diagram
\[
\xymatrix{
S^n\ar[r]^-{f_{1}}\ar[d]^-{c}			&A_1\ar[d]^-{c}\\
S^n\vee S^n\ar[r]^-{f_{1}\vee f_{1}}	&A_1\vee A_1
}
\]
homotopy commutes and fits in the middle of the following homotopy commutative diagram:
\[
\xymatrix{
S^n\ar[r]^-{c}\ar[d]^-{c}	&S^n\vee S^n\ar[r]^-{f_1\vee f_2}\ar[d]^-{c\vee 1}	&A_1\vee A_2\ar[d]^-{c\vee 1}\ar@/^1pc/[drr]^-{\Phi}	&	&\\
S^n\vee S^n\ar[r]^-{1\vee c}    &S^n\vee S^n\vee S^n\ar[r]^-{f_1\vee f_1\vee f_2}	&A_1\vee A_1\vee A_2\ar[r]^-{1\vee g\vee1}	&A_1\vee A_2\vee A_2\ar[r]^-{1\vee\triangledown}	&A_1\vee A_2.
}
\]
Here, the left-hand square homotopy commutes due to the coassociativity of the comultiplication $c$, and $\Phi$ is defined by the right-hand triangle.

Let $F$, $F'\colon S^n\rightarrow A_1\vee A_2$ be the maps defined by $F=(f_1\vee f_2)\circ c\simeq\jmath_1\circ f_1+\jmath_2\circ f_2$ and~$F'=\jmath_1\circ f_1+\jmath_2\circ (f_2+g\circ f_1)$. Since $F'$ is obtained by following the diagram anticlockwise, we obtain a homotopy $\Phi\circ F\simeq F'$. This witnesses the homotopy commutativity of the following diagram in which the rows are cofibre sequences and $\widetilde\Phi$ is an induced map of cofibres:
\[
\xymatrix{
S^n\ar[rr]^-{F}\ar@{=}[d]		&&A_1\vee A_2\ar[r]\ar[d]^-{\Phi}	&C(f_1,f_2)\ar[d]^-{\widetilde{\Phi}}\\
S^n\ar[rr]^-{F'}	&&A_1\vee A_2\ar[r]				&C(f_1,f_2+g\circ f_1).
}
\]
Since $\Phi$ is a homology equivalence, the five lemma implies that $\widetilde{\Phi}$ is as well. Since all spaces are simply connected, the Whitehead Theorem implies that $\tilde{\Phi}$ is a homotopy equivalence.
\end{proof}

\subsection{Some homotopy groups}

For $n\geq 3$ and $r,s\geq 1$ let $\imath_{rs}\colon P^n(2^r)\to P^n(2^s)$ be defined as follows. 
In case $r\geq s$, it is defined as the induced map of cofibres in the diagram
\begin{equation}\label{diagram_imath ell>k}
\begin{gathered}
\xymatrix{
S^{n-1}\ar[r]^-{2^r}\ar[d]^-{2^{r-s}}	&S^{n-1}\ar[r]^-{\imath_r}\ar@{=}[d]	&P^n(2^r)\ar[d]^-{\imath_{rs}}\\
S^{n-1}\ar[r]^-{2^s}					&S^{n-1}\ar[r]^-{\imath_s}&P^n(2^s).
}
\end{gathered}
\end{equation}
where $\imath_r\colon S^{n-1}\to P^n(2^r)$ is the inclusion of the bottom cell.
When $r<s$, the map $\imath_{rs}$ is similarly defined by the diagram
\begin{equation}\label{diagram_imath ell<k}
\begin{gathered}
\xymatrix{
S^{n-1}\ar[r]^-{2^r}\ar@{=}[d]	&S^{n-1}\ar[r]^-{\imath_r}\ar[d]^-{2^{s-r}}	&P^n(2^r)\ar[d]^-{\imath_{rs}}\\
S^{n-1}\ar[r]^-{2^s}&S^{n-1}\ar[r]^-{\imath_s}&P^n(2^s).
}
\end{gathered}
\end{equation}

\begin{lemma}\label{lemma_hmtpy gps Moore space}
Let $\eta\colon S^{n+1}\to S^n$ be the Hopf map.
\begin{enumerate}[label=\normalfont(\roman*)]
\item\label{hmtpy gp_pi 4 odd Moore sp}
If $t$ is odd, then $\pi_4(P^4(t))=\pi_5(P^4(t))=0$.
\item\label{hmtpy gp_pi 4 mod 2 Moore sp}
If $t=2^r$, then $\pi_4(P^4(2^r))\cong\Z/2\langle{\imath_r\circ\eta}\rangle$ and
\[
\pi_4(P^3(2^r))\cong\begin{cases}
\Z/4\langle{\xi_1}\rangle			&\text{for }r=1\\
\Z/2\langle{\xi_r}\rangle\oplus\Z/2\langle{\imath_r\circ\eta^2}\rangle	&\text{for }r\geq2
\end{cases}
\]
where $\xi_r$ is a lift of $\eta\colon S^4\to S^3$ through $P^3(2^r)$ and satisfies
\[
2\xi_1\simeq\imath_1\circ\eta^2
\quad\text{and}\quad
\xi_r=\imath_{1r}\circ\xi_1.
\]
\item\label{hmtpy gp Moore sp suspension}
The suspension $\Sigma\colon\pi_n(P^{n-1}(2^r))\to\pi_{n+1}(P^{n}(2^r))$ is an isomorphism for $n\geq 4$.
\end{enumerate}
\end{lemma}

\begin{proof}
The homotopy groups in~\ref{hmtpy gp_pi 4 odd Moore sp} and~\ref{hmtpy gp_pi 4 mod 2 Moore sp} are given in~\cite[Chapter 11.1 and Theorem~11.5.9]{baues}.

For $n\geq 5$, \ref{hmtpy gp Moore sp suspension} follows from the Freudenthal suspension theorem. We need to show that it holds for $n=4$. By~\cite[Proposition 11.1.12]{baues}, the double suspension
\[
\pi_4(P^3(2^r))\overset{\Sigma}{\longrightarrow}\pi_5(P^4(2^r))\overset{\Sigma}{\longrightarrow}\pi_6(P^5(2^r))
\]
is an isomorphism. Since the second $\Sigma$ is an isomorphism, so is the first.
\end{proof}

For $r\ge 1$, let $\C\PP^2(2^r)$ be the mapping cone of $S^3\overset{\eta}{\longrightarrow}S^2\overset{\imath_r}{\longrightarrow} P^3(2^r)$ and let {$\jmath_r\colon P^3(2^r)\to\C\PP^2(2^r)$ be the inclusion.

\begin{lemma}\label{lemma_pi_5 Sigma CP^2(2^r)}
Let $q\colon\Sigma\C\PP^2\to S^5$ and $q_r\colon \Sigma\C\PP^2(2^r)\to S^5$ be the quotient maps. Then
\begin{enumerate}[label=\normalfont(\roman*)]
\item\label{enumerate_pi_4 Sigma CP^2}
$\pi_4(\Sigma\C\PP^2)=\pi_4(\Sigma\C\PP^2(2^r))=0$;
\item\label{enumerate_pi_5 Sigma CP^2}
$\pi_5(\Sigma\C\PP^2)\cong\Z\langle{\alpha}\rangle$ and $\pi_5(\Sigma\C\PP^2(2^r))\cong\Z\langle{\alpha_r}\rangle\oplus\Z/2\langle{\epsilon_r}\rangle$,
where $\alpha$ and $\alpha_r$ satisfy
\[
q_*(\alpha)=(q_r)_*(\alpha_r)= 2\in\pi_5(S^5),
\]
and $\epsilon_r=\Sigma\jmath_r\circ\xi_r$;
\item\label{enumerate_map jmath Sigma CP^2(2) to CP^2(2)}
for $r,s\ge 1$ there exists a map $\jmath_{rs}\colon\Sigma\C\PP^2(2^r)\to\Sigma\C\PP^2(2^s)$ such that $\jmath_{rs}\circ\Sigma\jmath_r\simeq \Sigma\jmath_s\circ\imath_{rs}$ and 
\[
\jmath_{rs}\circ\epsilon_r\simeq
\begin{cases}
\ast        &r>s\\
\epsilon_s  &r\leq s.
\end{cases}
\]
\end{enumerate}
\end{lemma}

\begin{proof}
We first show~\ref{enumerate_pi_4 Sigma CP^2} and~\ref{enumerate_pi_5 Sigma CP^2}.
Applying the Blakers--Massey theorem to the cofibration
\[
S^4\overset{\imath_r\circ\eta}{\longrightarrow}P^4(2^r)\overset{\Sigma\jmath_r}{\longrightarrow}\Sigma\C\PP^2(2^r)
\]
yields an exact sequence
\begin{equation}\label{eqn_LES to compute pi Sigma CP(2^r)}
\xymatrix{
\pi_5(S^4)\ar[r]^-{(\imath_r\circ\eta)_*}	&\pi_5(P^4(2^r))\ar[r]^-{(\Sigma\jmath_r)_*}		&\pi_5(\Sigma\C\PP^2(2^r))\ar[r]^-{\partial}	&\pi_4(S^4)	&\\
	&\hspace{1cm} \ar[r]^-{(\imath_r\circ\eta)_*}	&\pi_4(P^4(2^r))\ar[r]^-{(\Sigma\jmath_r)_*}	&\pi_4(\Sigma\C\PP^2(2^r))\ar[r]	&0.
}
\end{equation}
By Lemma~\ref{lemma_hmtpy gps Moore space}, $\pi_4(P^4(2^r))\cong\Z/2\langle\imath_r\circ\eta\rangle$ so the second $(\imath_r\circ\eta)_*$ is surjective. It follows that $\pi_4(\Sigma\C\PP^2(2^r))$ is trivial and $\pi_5(\Sigma\C\PP^2(2^r))$ contains a $\mathbb{Z}$ summand with a generator $\alpha_r$ satisfying $\partial(\alpha_r)=2\in\pi_4(S^4)$. Since the boundary map $\partial$ is induced by the quotient map~$q_r$, there is a homotopy $q_r\circ\alpha_r\simeq 2\colon S^5\to S^5$.

We extract from \eqref{eqn_LES to compute pi Sigma CP(2^r)} the short exact sequence
\[
0\longrightarrow coker(\imath_r\circ\eta)_*\overset{(\Sigma\jmath_r)_*}{\longrightarrow}\pi_5(\Sigma\C\PP^2(2^r))\overset{\partial}{\longrightarrow}Im(\partial)\cong\Z\langle\alpha_r\rangle\longrightarrow0
\]
where $(\imath_r\circ\eta)_*$ is the first map appearing in~\eqref{eqn_LES to compute pi Sigma CP(2^r)}.
We will show that $coker(\imath_r\circ\eta)_*\cong\Z/2$, which then implies that $\pi_5(\Sigma\C\PP^2(2^r))\cong\Z\oplus\Z/2$.

In~\eqref{eqn_LES to compute pi Sigma CP(2^r)}, the first $(\imath_r\circ\eta)_*$ takes the generator $\eta$ of $\pi_5(S^4)\cong\Z/2$ to $\imath_r\circ\eta^2\in \pi_5(P^4(2^r))$. By Lemma~\ref{lemma_hmtpy gps Moore space}
\[
\pi_5(P^4(2^r))\cong\begin{cases}
\Z/4\langle{\xi_1}\rangle												&\text{if }r=1\\
\Z/2\langle{\xi_r}\rangle\oplus\Z/2\langle{\imath_r\circ\eta^2}\rangle	&\text{if }r\geq2.
\end{cases}
\]
When $r\geq2$, we immediately obtain $coker(\imath_r\circ\eta)_*\cong\Z/2\langle\xi_r\rangle$. When $r=1$, as $\imath_1\circ\eta^2\simeq2\xi_1$, $coker(\imath_1\circ\eta)_*$ is generated by the quotient image of $\xi_1$, which has order two. Let $\epsilon_r=\Sigma\jmath_r\circ\xi_r$ for $r\geq1$. Then $\epsilon_r\in\pi_5(\Sigma\C\PP^2(2^r))$ has order two and we have
\[
\pi_5(\Sigma\C\PP^2(2^r))\cong\Z\langle\alpha_r\rangle\oplus\Z/2\langle\epsilon_r\rangle.
\]

Similarly, one applies the Blakers--Massey theorem to the cofibration
\[
S^4\overset{\eta}{\longrightarrow}S^3\overset{\Sigma\jmath}{\longrightarrow}\Sigma\C\PP^2
\]
and see that $\pi_4(\Sigma\C\PP^2)$ is trivial and $\pi_5(\Sigma\C\PP^2)\cong\Z\langle\alpha\rangle$, where the generator $\alpha$ satisfies $q\circ\alpha\simeq 2\colon S^5\rightarrow S^5$.

Next we show~\ref{enumerate_map jmath Sigma CP^2(2) to CP^2(2)}. Given $r,s\in\N$, consider the diagram
\[
\xymatrix{
S^4\ar[r]^-{\eta}\ar[d]^-{2^{t}} &S^3\ar[r]^-{\imath_r}\ar[d]^-{2^{t}}    &P^4(2^r)\ar[d]^-{\imath_{rs}}\\
S^4\ar[r]^-{\eta}       &S^3\ar[r]^-{\imath_s}          &P^4(2^s)
}
\]
where $t=0$ for $r\geq s$ and $t=s-r$ for $r<s$.
The left square homotopy commutes since $\eta$ here is a suspension, and the right square homotopy commutes due to~\eqref{diagram_imath ell>k} and~\eqref{diagram_imath ell<k}.
Extend the diagram so as to obtain
\begin{equation}\label{diagram_Sigma CP to Sigma CP}
\begin{gathered}
\xymatrix{
S^4\ar[r]^-{\imath_r\circ\eta}\ar[d]^-{2^{t}}	&P^4(2^r)\ar[r]^-{\Sigma\jmath_r}\ar[d]^-{\imath_{rs}}	&\Sigma\C\PP^2(2^r)\ar[d]^-{\jmath_{rs}}\\
S^4\ar[r]^-{\imath_s\circ\eta}		&P^4(2^s)\ar[r]^-{\Sigma\jmath_s}			&\Sigma\C\PP^2(2^s)
}
\end{gathered}
\end{equation}
where the rows are cofibre sequences, and $\jmath_{rs}$ is an induced map. We show that $\jmath_{rs}$ has the asserted property. 

There is a string of homotopies
\[
\jmath_{rs}\circ\epsilon_r
\simeq\jmath_{rs}\circ(\Sigma\jmath_r\circ\xi_r)
\simeq\jmath_{rs}\circ\Sigma\jmath_r\circ(\imath_{1r}\circ\xi_1)
\simeq(\Sigma\jmath_{s}\circ\imath_{rs})\circ\imath_{1r}\circ\xi_1,
\]
where the first two homotopies are due to the definitions of $\epsilon_r$ and $\xi_r$, and the last is due to the right square of~\eqref{diagram_Sigma CP to Sigma CP}.
Since
\[
\imath_{rs}\circ\imath_{1r}\cong
\begin{cases}
2^{r-s}\imath_{1s}  &\text{if }r> s\\
\imath_{1s}         &\text{if }r\leq s,
\end{cases}
\]
the composite $\jmath_{rs}\circ\epsilon_r$ is null homotopic if $r>s$, and is homotopic to $\epsilon_s$ if $r\leq s$.
\end{proof}

\begin{remark}
The torsion generator $\epsilon_r\in\pi_5(\Sigma\C\PP^2(2^r))$ is a suspension since $\Sigma\jmath_r$ and $\xi_r$ are both suspensions by Lemma~\ref{lemma_hmtpy gps Moore space}.
\end{remark}
 
\subsection{A test for the nullity of Whitehead products}

For spaces $A$ and $B$ let
\[
[1,1]\colon\Sigma A\wedge B\rightarrow \Sigma A\vee \Sigma B
\]
denote the external Whitehead product whose cofibre is $\Sigma A\times\Sigma B$.

\begin{lemma}\label{lemma_WH nullity criterion}
For $n\geq 2$ let $f\colon S^n\to\Sigma A\wedge B$ be a map, and let $C$ be the mapping cone of the composite
\[
S^n\overset{f}{\longrightarrow}\Sigma A\wedge B\overset{[1,1]}{\longrightarrow}\Sigma A\vee\Sigma B.
\]
Suppose $\Sigma A\wedge B$ is $(n-1)$-connected and $H_n(\Sigma A\wedge B;\Z)$ is a cyclic group. If all cup products in $H^*(C;R)$ are trivial for any ring $R$, then $f\simeq*$.
\end{lemma}

\begin{proof}
Since $\Sigma A\wedge B$ is $(n-1)$-connected, the Hurewicz theorem implies that $f$ is null homotopic if and only if the induced map
\[
f_*\colon H_n(S^n;\Z)\to H_n(\Sigma A\wedge B;\Z)
\]
is trivial. By assumption $H_n(\Sigma A\wedge B;\Z)\cong\Z/k$ for $k\in\{2,3,\ldots,\infty\}$, where $\Z/{\infty}$ means $\Z$. The universal coefficient theorem then implies that $f_*$ is trivial if and only if
\[
f^*\colon H^n(\Sigma A\wedge B;\Z/k)\to H^n(S^n;\Z/k)
\]
is trivial. We show that the latter holds under the assumption.

Consider the commutative diagram
\[
\xymatrix{
S^n\ar[r]^-{f}\ar@{=}[d]			&\Sigma A\wedge B\ar[r]^-{q}\ar[d]^-{[1,1]}	&C_f\ar[d]\\
S^n\ar[r]^-{[1,1]\circ f}\ar[d]	&\Sigma A\vee\Sigma B\ar[r]\ar[d]	&C\ar[d]^-{\delta}\\
\ast\ar[r]							&\Sigma A\times\Sigma B\ar@{=}[r]	&\Sigma A\times\Sigma B
}
\]
where $q$ is the quotient map, $\delta$ is an induced map, and all columns and rows are cofibre sequences. Apply $H^*(-;\Z/k)$ and suppress coefficients to obtain a commutative diagram with exact rows and columns
\[
\xymatrix{
0\ar[r]	&H^n(C_f)\ar[r]^-{q^*}\ar[d]^-{\epsilon}	&H^n(\Sigma A\wedge B)\ar[r]^-{f^*}\ar[d]^-{\epsilon'}	&H^n(S^n)\\
	&H^{n+1}(\Sigma A\times\Sigma B)\ar@{=}[r]\ar[d]^-{\delta^*}	&H^{n+1}(\Sigma A\times\Sigma B)\ar[d]	&\\
	&H^{n+1}(C)\ar[r]	&H^{n+1}(\Sigma A\vee\Sigma B)	&
}
\]
where $\epsilon$ and $\epsilon'$ are connecting maps. Let $\alpha\in H^n(\Sigma A\wedge B;\Z/k)\cong\Z/k$ be a generator. To show the triviality of $f$ we need to show that $\alpha$ is in $\im(q^*)$.

Notice that $\epsilon'$ is injective and maps $H^n(\Sigma A\wedge B)\cong H^{n+1}(\Sigma A\wedge\Sigma B)$ onto the subgroup of~$H^{n+1}(\Sigma A\times\Sigma B)$ generated by cup products of the form $pr_1^*(x)\cup pr_2^*(y)$, where $pr_1\colon\Sigma A\times\Sigma B\to\Sigma A$ and $pr_2\colon\Sigma A\times\Sigma B\to\Sigma B$ are projections. Since all cup products in $H^*(C)$ are assumed to be trivial, we therefore have $\delta^*(\epsilon'(\alpha))=0$. Now $\epsilon=\epsilon'\circ q^*$ is injective, so there is a unique class $\widetilde\alpha\in H^{n+1}(C_f)$ such that $\epsilon(\widetilde\alpha)=\epsilon'(\alpha)$. In particular $\alpha=q^*(\widetilde\alpha)$.
\end{proof}

\section{Strategy for proving the main theorem}\label{sect3}
Suppose the homology groups of $X$ are given by~\eqref{table_original M hmlgy}.
Let $\varphi\colon\bigvee^m_{i=1}S^1\to M$ be a map inducing an isomorphism on $H_1(\;)$, and let $X'$ be the mapping cone of $\varphi$. Then the homology groups of $X'$ are given by
\[
\begin{tabular}{C{1.4cm}|C{1.4cm}|C{1.4cm}|C{1.4cm}|C{1.4cm}|C{1.4cm}|C{1.4cm}}
$i$	&$0$	&$1$	&$2$	&$3$	&$4$	&$5$\\
\hline
$H_i(X')$	&$\Z$	&$0$	&$\Z^n\oplus T$	&$\Z^{n}$	&$\Z^m$	&$\Z$
\end{tabular}.
\]
Using the argument in \cite[Lemma 5.1]{ST19} we obtain a homotopy equivalence
\begin{equation}\label{equation_Sigma X'}
\Sigma X\simeq\bigvee^m_{i=1}S^2\vee\Sigma X'.
\end{equation}
Hence it suffices to study the homotopy type of $\Sigma X'$.
Since $\Sigma X'$ is simply connected, it has a minimal cell structure which can be constructed as follows.
\begin{itemize}
\item
Let $W_2=\bigvee^n_{i=1}S^3\vee\bigvee^{\ell}_{j=1}P^4(t_j)$;

\item
For $3\leq i\leq 5$, let $W_i$ be the mapping cone of a map
\begin{equation}\label{eqn_attaching map f_5}
\varphi_i\colon\bigvee^{n_i}_{j=1}S^i\longrightarrow W_{i-1}
\end{equation}
where $n_i=\text{rank}(H_i(X))$ is the $i$-th Betti number of $X$ and $\varphi_i$ induces a trivial homomorphism in homology;

\item 
Then $W_5\simeq\Sigma X'$ by \cite[Theorem 7.3.2]{arkowitz}. 
\end{itemize}

We will show that each of the $W_i$ is homotopy equivalent to a wedge of recognisable spaces, and hence obtain Theorem~\ref{thm_main theorem}. First we begin with $W_3$.

\begin{lemma}\label{lemma_W_3}
There is a homotopy equivalence
\[
W_3\simeq\bigvee^n_{i=1}(S^3\vee S^4)\vee\bigvee^{\ell}_{j=1}P^4(t_j).
\]
\end{lemma}

\begin{proof}
Since $W_3$ is the mapping cone of $\varphi_3\colon\bigvee^n_{i=1}S^3\to W_2$ given in~\eqref{eqn_attaching map f_5}, it suffices to show that $\varphi_3$ is null homotopic. Consider the commutative diagram
\[
\xymatrix{
\pi_3(\bigvee^n_{i=1}S^3)\ar[r]^-{[\varphi_3]_*}\ar[d]^-{h}	&\pi_3(W_2)\ar[d]^-{h}\\
H_3(\bigvee^n_{i=1}S^3)\ar[r]^-{(\varphi_3)_*}					&H_3(W_2)
}
\]
where the vertical maps labelled $h$ are Hurewicz morphisms. Since $W_2$ is 2-connected, the Hurewicz Theorem implies that the two $h$'s are isomorphisms. Since $(\varphi_3)_*$ is trivial by assumption, $\varphi_3\simeq[\varphi_3]_*(id)$ is null homotopic.
\end{proof}

The calculations needed for $W_4$ and $W_5$ are significantly longer and are completed in the next section.

\section{Homotopy types of $W_4$ and $W_5$}\label{section_homotopy types of W_4 and W_5}

Let $\mathscr{W}$ be the collection of spaces that are of the form
\[
(S^3)^{\vee n_3}\vee(S^4)^{\vee n_4}\vee(S^5)^{\vee n_5}\vee\bigvee^{a}_{i=1}P^4(2^{r_i})\vee\bigvee_{t\in \mathcal{T}}P^4(t)\vee(\Sigma\C\PP^2)^{\vee b}\vee\bigvee^{c}_{j=1}\Sigma\C\PP^2(2^{s_j})
\]
where $n_3,n_4,n_5,b$ are non-negative integers, $\mathcal{T}$ is a collection of odd numbers, and $\{r_i\}^a_{i=1}$ and $\{s_j\}^c_{j=1}$ are sequences of positive integers such that $r_i\leq r_{i+1}$ and $s_j\leq s_{j+1}$.
We will show that $W_4$ is homotopy equivalent to a wedge in $\mathscr{W}$ in Lemma~\ref{lemma_homotopy type of W_4}, and compute the homotopy type of $W_5$ in Lemma~\ref{lemma_homotopy type of C_varphi for a=5}.

We begin with a useful lemma for later calculations.

\begin{lemma}\label{lemma_varphi simplify no WH prod}
Fix a wedge $W\in\mathscr{W}$. Let $\varphi\colon S^d\to W$ be a map where $d=4$ or $5$, and let~$C_{\varphi}$ be its mapping cone.
If $\varphi$ induces the trivial morphism in homology and all cup products in~$\tilde{H}^*(C_{\varphi};R)$ are trivial for any ring $R$, then there is a homotopy
\[
\varphi\simeq\sum_{A}\varphi(A)
\]
where $A$ runs over all wedge summands in $W$ such that $A\neq S^d$, 
and each $\varphi(A)$ is
\[
\varphi(A)\colon S^d\overset{\varphi}{\longrightarrow}W\overset{\text{pinch}}{\longrightarrow}A\hookrightarrow W.
\]
\end{lemma}

\begin{proof}
Let $\mathscr{A}$ be the collection of wedge summands in $W$.
Since each wedge summand $A$ is a suspension and is at least 2-connected, the Hilton--Milnor theorem implies that
\[
\varphi\simeq\sum_{A\in\mathscr{A}}\varphi(A)+
\sum_{\substack{A,A'\in\mathscr{A}\\ A\neq A'}}w(A,A')
\]
where $w(A,A')\colon S^d\to\Sigma^{-1}A\vee A'$ is the composite
\[
w(A,A')\colon S^d\overset{f}{\longrightarrow}\Sigma^{-1}A\wedge A'\overset{[1,1]}{\longrightarrow}A\vee A'\hookrightarrow W
\]
of a map $f\colon S^d\to\Sigma^{-1}A\wedge A'$, the Whitehead product $[1,1]$, and the inclusion $A\vee A'\hookrightarrow W$.

For those $A=S^d$, the map $\varphi(A)$ is null homotopic since $\varphi_*$ is trivial in homology by assumption. It remains to show that each $w(A,A')$ is null homotopic.

Suppose $d=4$. Note that each $\Sigma^{-1}A\wedge A'$ is at least 4-connected. Hence $f\colon S^4\to\Sigma^{-1}A\wedge A'$ is null homotopic and so is $w(A,A')$.

Suppose $d=5$. 
Consider the composite
$\varphi_{A,A'}\colon S^5\overset{\varphi}{\longrightarrow}W\overset{\text{pinch}}{\longrightarrow}A\vee A'$.
We have
\[
\varphi_{A,A'}\simeq\pi\circ\varphi(A)+\pi\circ\varphi(A')+[1,1]\circ f,
\]
where $\pi\colon W\to A\vee A'$ is the pinch map.
Let $C_A,C_{A'},C_{AA'}$ and $C_w$ be the mapping cones of~\mbox{$\pi\circ\varphi(A),\pi\circ\varphi(A'),\varphi_{A,A'}$} and $[1,1]\circ f$ respectively. Since cup products in $\tilde{H}^*(C_{\varphi};R)$ are trivial, so are the cup products in $\tilde{H}^*(C_{AA'};R)$. Note that $\pi\circ\varphi(A)$ and $\pi\circ\varphi(A')$ are suspensions, implying that cup products in $\tilde{H}^*(C_A;R)$ and $\tilde{H}^*(C_{A'};R)$ are trivial. Then cup products in~$\tilde{H}^*(C_w;R)$ are also trivial by \cite[Lemma 2.5]{CS} (the dimension assumption of the lemma is not necessary).
Since $\Sigma^{-1}A\wedge A'$ satisfies the conditions of Lemma~\ref{lemma_WH nullity criterion}, $f$ is null homotopic and so is $w(A,A')$.
\end{proof}

\subsection{The homotopy type of $W_4$}
Since $W_4$ is the mapping cone of
$\varphi_4\colon\bigvee^m_{i=1}S^4_i\to W_3$
given in~\eqref{eqn_attaching map f_5}, it can be obtained by adding a 5-cell to $W_3$ iteratively. Observe that
\[
W_3\simeq\bigvee^n_{j=1}(S^3\vee S^4)\vee\bigvee^{\ell}_{k=1}P^4(t_k)
\]
is a wedge in $\mathscr{W}$. To prove that $W_4$ is in $\mathscr{W}$, it suffices to show that the complex obtained by attaching a 5-cell is always homotopy equivalent to a wedge in $\mathscr{W}$.

Consider a wedge $W\in\mathscr{W}$ and label its 3-spheres, that is
\begin{equation}\label{eqn_wedge form}
W=\bigvee^{n_3}_{i=1}S^3_i\vee\bigvee^{n_4}_{j=1}S^4\vee\bigvee^{n_5}_{k=1}S^5\vee\bigvee^{a}_{u=1}P^4(2^{r_u})\vee\bigvee_{t\in \mathcal{T}}P^4(t)\vee\bigvee^{b}_{l=1}\Sigma\C\PP^2\vee\bigvee^{c}_{v=1}\Sigma\C\PP^2(2^{s_v}).
\end{equation}
Given any map $\varphi\colon S^4\to W$, we claim that its mapping cone $C_{\varphi}$ is also a wedge in $\mathscr{W}$ if $\varphi$ induces the trivial morphism in homology.
 
\begin{lemma}\label{lemma_varphi_4 simplify}
Let $W$ be a wedge in $\mathscr{W}$ as in~\eqref{eqn_wedge form} and let $\varphi\colon S^4\to W$ be a map.
If $\varphi$ induces the trivial morphism in homology, then
\[
\varphi\simeq\sum^{n_3}_{i=1}\varphi(S^3_i)+\sum^{a}_{u=1}\varphi(P^4(2^{r_u})).
\]
\end{lemma}

\begin{proof}
Notice that  all cup products in $\tilde{H}^*(C_{\varphi})$ are trivial as $H^i(C_{\varphi})=0$ for $i=1,2$.
By Lemma~\ref{lemma_varphi simplify no WH prod} it suffices to show that the components
\[
\varphi(S^5),\quad
\varphi(P^4(t)),\quad
\varphi(\Sigma\C\PP^2),\quad
\text{and}\quad
\varphi(\Sigma\C\PP^2(2^{s}))
\]
are null homotopic.
This holds due to the fact that $\pi_4(S^5)=0$ and Lemmas~\ref{lemma_hmtpy gps Moore space} and~\ref{lemma_pi_5 Sigma CP^2(2^r)}.
\end{proof}

Recall from Lemma~\ref{lemma_hmtpy gps Moore space} that
\[
\pi_4(S^3)\cong\Z/2\langle\eta\rangle
\quad\text{and}\quad
\pi_4(P^4(2^r))\cong\Z/2\langle\imath_r\circ\eta\rangle,
\]
where $\eta$ is the suspended Hopf map and $\imath_{r}\colon S^3\to P^4(2^r)$ is the inclusion of the bottom cell.
For $1\leq i\leq n_3$ and $1\leq u\leq a$, let $\eta_i$ and $\bar{\eta}_u$ be the composites
\[
\eta_i\colon S^4\overset{\eta}{\longrightarrow}S^3_i\hookrightarrow W
\quad\text{and}\quad
\bar{\eta}_u\colon S^4\overset{\imath_{r_u}\circ\eta}{\longrightarrow}P^4(2^{r_u})\hookrightarrow W.
\]
Then Lemma~\ref{lemma_varphi_4 simplify} implies that there are some $\Z/2$-coefficients $x_i,y_u$ such that
\begin{equation}\label{eqn_varphi_4 simplify}
\varphi\simeq\sum^{n_3}_{i=1}x_i\eta_i+\sum^{a}_{u=1}y_u\bar{\eta}_u.
\end{equation}
Following the idea in~\cite{CS, ST19}, we may choose a different map $\varphi'\colon S^4\to W$ such that $C_{\varphi}\simeq C_{\varphi'}$.

\begin{lemma}\label{lemma_htpy equiv of mapping cone of varphi_4}
Let $C(x_1,\ldots,x_{n_3};y_1,\ldots,y_{a})$ be the mapping cone of $\sum^{n_3}_{i=1}x_i\eta_i+\sum^{a}_{u=1}y_u\bar{\eta}_u$.
Then $C(x_1,\ldots,x_{n_3};y_1,\ldots,y_{a})$ is homotopy equivalent to
\begin{enumerate}[label=\normalfont(\roman*)]
\item\label{W_4 simply_ permute}
$C(x_{\sigma(1)},\ldots, x_{\sigma(n_3)};y_1,\ldots,y_{a})$ for any permutation $\sigma\colon\{1,\ldots,n_3\}\to\{1,\ldots,n_3\}$;
\item\label{W_4 simply_ x->x}
$C(x_1,\ldots, x_{l-1},x_k+x_l,x_{l+1},\ldots,x_{n_3};y_1,\ldots,y_{a})$ for some distinct $k,l\in\{1,\ldots,n_3\}$;
\item\label{W_4 simply_ x->y}
$C(x_1,\ldots, ,x_{n_3};y_1,\ldots,y_{l-1},x_k+y_l,y_{l+1},\ldots,y_{a})$ for some $k\in\{1,\ldots,n_3\}$ and $l\in\{1,\ldots,a\}$;
\item\label{W_4 simply_ y->y}
$C(x_1,\ldots, ,x_{n_3};y_1,\ldots,y_{l-1},y_k+y_l,y_{l+1},\ldots,y_{a})$ for some distinct $k,l\in\{1,\ldots,a\}$ such that $k\geq l$.
\end{enumerate}
\end{lemma}

\begin{proof}
First we prove~\ref{W_4 simply_ permute}. Let $\varphi'=\sum^{n_3}_{i=1}x_{\sigma(i)}\eta_i+\sum^{a}_{u=1}y_u\bar{\eta}_u$, and let $\sigma\colon W\to W$ be the map which permutes the wedge summands by mapping each $S^3_i$ to $S^3_{\sigma(i)}$ and mapping the other wedge summands onto themselves by the identity map. Then there is a commutative diagram of cofibre sequences
\[
\xymatrix{
S^4\ar[r]^-{\varphi}\ar@{=}[d]	&W\ar[r]\ar[d]^-{\sigma}	&C(x_1,\ldots,x_{n_3};y_1,\ldots,y_{a})\ar[d]^-{\tilde{\sigma}}\\
S^4\ar[r]^-{\varphi'}			&W\ar[r]			&C(x_{\sigma(1)},\ldots,x_{\sigma(n_3)};y_1,\ldots,y_{a})
}
\]
where $\tilde{\sigma}$ is an induced map. Since $\sigma$ is a homology isomorphism, so is $\tilde{\sigma}$ by the five lemma. Consequently, the Whitehead Theorem implies that $\tilde{\sigma}$ is a homotopy equivalence.

Next we prove~\ref{W_4 simply_ x->x} to~\ref{W_4 simply_ y->y}.
Use Lemma~\ref{lemma_simplify mapping cone} and take $g\colon A_k\to A_l$ to be
\begin{itemize}
\item
the identity map $A_k\to A_l$ if $A_k=A_l$
\item
the inclusion $\imath_{r_l}\colon S^3\to P^4(2^{r_l})$ if $A_k=S^3$ and $A_l=P^4(2^{r_l})$
\item
the map $\imath_{r_kr_l}\colon P^4(2^{r_k})\to P^4(2^{r_l})$ in~\eqref{diagram_imath ell>k} if $A_k=P^4(2^{r_k}),A_l=P^4(2^{r_l})$ and $k\geq l$.
\end{itemize}
Then the homotopy equivalences follow immediately.
\end{proof}

\begin{lemma}\label{lemma_hmtpy type of varphi for d=4}
Under the assumption of Lemma~\ref{lemma_varphi_4 simplify},
$C_{\varphi}$ is a wedge in $\mathscr{W}$.
\end{lemma}

\begin{proof}
By~\eqref{eqn_varphi_4 simplify} we have $C_{\varphi}\simeq (x_1,\ldots,x_{n_3};y_1,\ldots,y_{a})$ for some $\Z/2$ coefficients $x_i$ and $y_j$. We divide the proof into 3 cases.

\textit{Case 1}: If all $x_i$'s and $y_j$'s are zero, then $\varphi$ is null homotopic and $C_{\varphi}\simeq W\vee S^5$ is in $\mathscr{W}$;

\textit{Case 2}: Suppose $x_i=1$ for some $i\in\{1,\ldots,n_3\}$. By Lemma~\ref{lemma_htpy equiv of mapping cone of varphi_4}~\ref{W_4 simply_ permute} we assume it is $x_1$. If there are non-zero coefficients $x_j$ and $y_u$ in Equation~\eqref{eqn_varphi_4 simplify}, then Lemma~\ref{lemma_htpy equiv of mapping cone of varphi_4}~\ref{W_4 simply_ x->x} and~\ref{W_4 simply_ x->y} implies that $C(x_1,\ldots,x_{n_3};y_1,\ldots,y_{a})$ is homotopy equivalent to
\[
C(x_1,\ldots,x_{j-1},x_j+x_1,x_{j+1},\ldots,x_{n_3};y_1,\ldots,y_{u-1},y_u+x_1,y_{u+1},\ldots,y_{a}).
\]
Note that $x_j+x_1$ and $y_u+x_1$ are both zero. Repeat the argument to annihilate other non-zero coefficients. In the end we have $C_{\varphi}\simeq C(1,0,\ldots,0;0,\ldots,0)$. Observe that it is the mapping cone of $\eta_1$ and is homotopy equivalent to
\[
(S^3)^{\vee(n_3-1)}\vee(S^4)^{\vee n_4}\vee(S^5)^{\vee n_5}
\vee\bigvee^{a}_{u=1}P^4(2^{r_u})\vee\bigvee_{t\in \mathcal{T}}P^4(t)\vee(\Sigma\C\PP^2)^{\vee(b+1)}\vee\bigvee^{c}_{v=1}\Sigma\C\PP^2(2^{s_v}).
\]
Therefore $C_{\varphi}$ is in $\mathscr{W}$;

\textit{Case 3}: Suppose $x_i=0$ for all $i$ and $y_u=1$ for some $u\in\{1,\ldots,a\}$. Let $\mu$ be the largest index of those non-zero $y_u$. Then $C_{\varphi}=C(0,\ldots,0;y_1,\ldots,y_{\mu},0,\ldots,0)$
where $y_{\mu}=1$. Use Lemma~\ref{lemma_htpy equiv of mapping cone of varphi_4}~\ref{W_4 simply_ y->y} to annihilate other non-zero $y_u$. In the end we have
\[
C_{\varphi}\simeq C(0,\ldots,0;0,\ldots,0,y_{\mu},0,\ldots,0),
\]
which is the mapping cone of $\tilde{\eta}_{\mu}$ and is homotopy equivalent to
\[
(S^3)^{\vee n_3}\vee(S^4)^{\vee n_4}\vee(S^5)^{\vee n_5}
\vee\bigvee_{\substack{1\leq u\leq a \\ u\neq\mu}}P^4(2^{r_u})\vee\bigvee_{t\in \mathcal{T}}P^4(t)\vee(\Sigma\C\PP^2)^{\vee b}\vee\bigvee^{c}_{v=1}\Sigma\C\PP^2(2^{s_v})\vee\Sigma\C\PP^2(2^{r_{\mu}}).
\]
Therefore $C_{\varphi}$ is in $\mathscr{W}$.
\end{proof}

Now we have all the ingredients to show that $W_4$ is in $\mathscr{W}$.

\begin{lemma}\label{lemma_homotopy type of W_4}
There is a homotopy equivalence
\[
W_4\simeq
(S^3)^{\vee n_3}\vee(S^4)^{\vee n_4}\vee(S^5)^{\vee n_5}\vee\bigvee^{a}_{u=1}P^4(2^{r_u})\vee\bigvee_{t\in \mathcal{T}}P^4(t)\vee(\Sigma\C\PP^2)^{\vee b}\vee\bigvee^{c}_{v=1}\Sigma\C\PP^2(2^{s_v})
\]
for some integers $n_3,n_4,n_5,a,b,c,r_u,s_v$ and some collection $\mathcal{T}$ of odd numbers.
\end{lemma}

\begin{proof}
Let $\{Y_i\}_{1\leq i\leq m}$ be a sequence of spaces where $Y_0=W_3$ and $Y_i$ is the mapping cone of
\[
f_i\colon S^4_i\hookrightarrow\bigvee^m_{j=1} S^4_j\overset{\varphi_4}{\longrightarrow}W_3\hookrightarrow Y_{i-1}.
\]
for $1\leq i\leq m$.
Then $Y_m\simeq W_4$.
We prove that each $Y_i$ is in $\mathscr{W}$ by induction on $i$.

By Lemma~\ref{lemma_W_3}, $Y_0=W_3$ is homotopy equivalent to
\[
W_3\simeq\bigvee^n_{j=1}(S^3\vee S^4)\vee\bigvee_{t\in T} P^4(t),
\]
where $T=\{p^r\}$ is a collections of powers of primes.
Let
\[
T_{odd}=\{p^r\in T\mid \text{$p$ is an odd prime}\}
\quad\text{and}\quad
T_{even}=\{p^r\in T\mid p=2\}
\]
be the subsets consisting of powers of odd and even primes respectively. Further, let
\[
T_{even}=\{2^{r_1},\ldots,2^{r_{a}}\}
\]
such that $r_u\leq r_{u+1}$. Therefore $Y_0\simeq\bigvee^{m}_{j=1}(S^3\vee S^4)\vee\bigvee^{a}_{u=1}P^4(2^{r_u})\vee\bigvee_{t\in T_{odd}}P^4(t)$
and hence is in $\mathscr{W}$.

Assume that $Y_i$ is in $\mathscr{W}$. Let $Y_{i+1}$ be the mapping cone of $f_{i+1}\colon S^4\to Y_i$.
Since $\varphi_4$ induces the trivial morphism in homology, so does $f_i$.
Lemmas~\ref{lemma_varphi_4 simplify} and~\ref{lemma_hmtpy type of varphi for d=4} imply that
$Y_{i+1}$ is in $\mathscr{W}$.

By induction all $Y_i$' s are in $\mathscr{W}$ and, in particular, so is $W_4=Y_m$.
\end{proof}

\subsection{The homotopy type of $W_5$}
Recall that $W_5$ is the mapping cone of
$\varphi_5\colon S^5\to W_4$
given in~\eqref{eqn_attaching map f_5} such that $\varphi_5$ induces the trivial morphism in homology. Since $W_5\simeq\Sigma X'$ is a suspension, all cup products in $\tilde{H}^*(W_5)$ are trivial.

\begin{lemma}\label{lemma_varphi_5 simplify}
Label the $3$-spheres and $4$-spheres in $W_4$ by $S^3_i$ and $S^4_j$. Then
\[
\varphi_5\simeq\sum^{n_3}_{i=1}\varphi(S^3_i)+\sum^{n_4}_{j=1}\varphi(S^4_j)+\sum^{a}_{u=1}\varphi(P^4(2^{r_u}))+\sum^{c}_{v=1}\varphi(\Sigma\C\PP^2(2^{s_v})).
\]
Furthermore, each $\varphi(\Sigma\C\PP^2(2^{s_v}))$ is in $\Z/2\langle\epsilon_{s_v}\rangle$.
\end{lemma}

\begin{proof}
Apply Lemma~\ref{lemma_varphi simplify no WH prod} to obtain $\varphi_5\simeq\sum_{A\neq S^5}\varphi(A)$.
By Lemma~\ref{lemma_hmtpy gps Moore space} each $\varphi(P^4(t))$ is null homotopic for $t$ odd, so we have
\[
\varphi_5\simeq\sum^{n_3}_{i=1}\varphi(S^3_i)+\sum^{n_4}_{j=1}\varphi(S^4_j)+\sum^{a}_{u=1}\varphi(P^4(2^{r_u}))+\sum^b_{k=1}\varphi(\Sigma\C\PP^2)+\sum^{c}_{v=1}\varphi(\Sigma\C\PP^2(2^{s_v})).
\]
By Lemma~\ref{lemma_pi_5 Sigma CP^2(2^r)} each $\varphi(\Sigma\C\PP^2)\simeq A\alpha$ and $\varphi(\Sigma\C\PP^2(2^{s_v}))\simeq B\alpha_{s_v}+C\epsilon_{s_v}$ for some $A,B\in\Z$ and $C\in\Z/2$. It remains to show that $A=B=0$.

Let $q_{s_v}\colon\Sigma\C\PP^2(2^{s_v})\to S^5$ be the quotient map. Observe that the composite
\[
S^5\overset{\varphi(\Sigma\C\PP^2(2^{s_v}))}{\longrightarrow}W_4\overset{\text{pinch}}{\longrightarrow}\Sigma\C\PP^2(2^{s_v})\overset{q_{s_v}}{\longrightarrow}S^5
\]
is homotopic to $q_{s_v}\circ\alpha_{s_v}$ which has degree 2 by Lemma~\ref{lemma_pi_5 Sigma CP^2(2^r)}. The induced morphism
\[
\varphi(\Sigma\C\PP^2(2^{s_v}))_*\colon H_5(S^5)\cong\Z\to H_5(\Sigma\C\PP^2(2^{s_v}))\cong\Z
\]
is given by $1\mapsto\pm 2B$. Since by assumption $\varphi_5$ induces the trivial homology morphism, so does $\varphi(\Sigma\C\PP^2(2^{s_v}))$ and $B$ has to be zero. A similar argument shows that $A=0$. The lemma then follows.
\end{proof}

Lemma~\ref{lemma_varphi_5 simplify} says that $\varphi_5$ decomposes into a sum of maps that factor through the wedge summands
$P^4(2^{r_u})$,
$S^3_i$,
$S^4_j$,
$\Sigma\C\PP^2(2^{s_v})$
of $W_4$. Define $\tilde{\eta}_i$, $\tilde{\xi}_j$, $\tilde{\epsilon}_k\colon S^5\to W_4$
as follows:
\begin{itemize}
\item
$\tilde{\eta}_i$ is
$\begin{cases}
S^5\overset{\imath\circ\eta^2}{\longrightarrow}P^4(2^{r_i})\hookrightarrow W_4	&\text{for }1\leq i\leq a\\
S^5\overset{\eta^2}{\longrightarrow}S^3_{i-a}\hookrightarrow W_4	&\text{for }a+1\leq i\leq a+n_3\\
S^5\overset{\eta}{\longrightarrow}S^4_{i-a-n_3}\hookrightarrow W_4	&\text{for }a+n_3+1\leq i\leq a+n_3+n_4
\end{cases}$
\item
$\tilde{\xi}_j$ is $S^5\overset{\xi_{r_j}}{\longrightarrow}P^4(2^{r_j})\hookrightarrow W_4$ for $1\leq j\leq a$
\item
$\tilde{\epsilon}_k$ is $S^5\overset{\epsilon_{s_k}}{\longrightarrow}\Sigma\C\PP^2(2^{s_k})\hookrightarrow W_4$ for $1\leq k\leq c$.
\end{itemize}
To keep track of which wedge summands these maps factor through, we define $\{r(i)\}^{1+n_3+n_4}_{i=1}$ and $\{s(j)\}^c_{j=1}$ by
\[
r(i)=\begin{cases}
r_i		&\text{for }1\leq i\leq a\\
\omega	&\text{for }a+1\leq i\leq a+n_3\\
\omega+1	&\text{for }a+n_3+1\leq i\leq a+n_3+n_4
\end{cases}
\quad\text{and}\quad
s(j)=s_j\quad\text{for }1\leq j\leq c,
\]
where $\omega$ is the first infinite ordinal.
By Lemmas~\ref{lemma_hmtpy gps Moore space} and~\ref{lemma_pi_5 Sigma CP^2(2^r)}, we have
\begin{equation}\label{eqn_eta xi epsilon restrictions}
2\tilde{\eta}_i\simeq\ast,\quad
2\tilde{\epsilon}_k\simeq\ast,\quad
2\tilde{\xi}_j\simeq\begin{cases}
\tilde{\eta}_j  &\text{if }r(j)=1\\
\ast            &\text{if }r(j)>1.
\end{cases}
\end{equation}
It follows that the attaching map $\varphi_5$ can be written as a linear combination
\begin{equation}\label{eqn_varphi_5 simplify}
\varphi_5\simeq\sum^{a+n_3+n_4}_{i=1}x_i\tilde{\eta}_i+\sum^{a}_{j=1}y_j\tilde{\xi}_j+\sum^{c}_{k=1}z_k\tilde{\epsilon}_k.
\end{equation}
for some coefficients $x_i,y_j,z_k\in\{0,1\}$. We now establish a lemma similar to Lemma~\ref{lemma_htpy equiv of mapping cone of varphi_4} for~$\varphi_5$.

\begin{lemma}\label{lemma_htpy equiv of mapping cone of varphi_5}
Let $\mathbf{x}=(x_i)_{1\leq i\leq a+n_3+n_4},
\mathbf{y}=(y_j)_{1\leq j\leq a}$, and $\mathbf{z}=(z_k)_{1\leq k\leq c}$ be sequences of non-negative integers,
and let $C(\mathbf{x},\mathbf{y},\mathbf{z})$ be the mapping cone of
\[
\sum^{a+n_3+n_4}_{i=1}x_i\tilde{\eta}_i+\sum^{a}_{j=1}y_j\tilde{\xi}_j+\sum^{c}_{k=1}z_k\tilde{\epsilon}_k.
\]
Then the homotopy type of $C(\mathbf{x},\mathbf{y},\mathbf{z})$ is restricted by the following conditions:
\begin{enumerate}[label=\normalfont(\roman*)]
\item\label{eta_i has order 2}
if $x_i=2$, then
$C(\mathbf{x},\mathbf{y},\mathbf{z})\simeq C(x_1,\ldots,x_{i-1},0,x_{i+1},\ldots,x_{a+n_3+n_4};\mathbf{y};\mathbf{z})$;
\item\label{xi_j has order 4 if r_j=1}
if $y_j=2$ with $r(j)=1$, then\\
$C(\mathbf{x},\mathbf{y},\mathbf{z})\simeq C(x_1,\ldots,x_{j-1},x_j+1,x_{j+1},\ldots,x_{a+n_3+n_4};y_1,\ldots,y_{j-1},0,y_{j+1},\ldots,y_{a};\mathbf{z})$;
\item\label{y_j=3} if $y_j=3$ with $r(j)=1$, then\\
$C(\mathbf{x},\mathbf{y},\mathbf{z})\simeq C(x_1,\ldots,x_{j-1},x_j+1,x_{j+1},\ldots,x_{a+n_3+n_4};y_1,\ldots,y_{j-1},1,y_{j+1},\ldots,y_{a};\mathbf{z})$;
\item\label{xi_j has order 2 if r_j>1}
if $y_j=2$ with $r(j)>1$, then
$C(\mathbf{x},\mathbf{y},\mathbf{z})\simeq C(\mathbf{x};y_1,\ldots,y_{j-1},0,y_{j+1},\ldots,y_{a};\mathbf{z})$;
\item\label{epsilon_k has order 2}
if $z_k=2$, then
$C(\mathbf{x},\mathbf{y},\mathbf{z})\simeq C(\mathbf{x};\mathbf{y}; z_1,\ldots,z_{k-1},0,y_{k+1},\ldots,z_{c})$.
\end{enumerate}
Furthermore, $C(\mathbf{x},\mathbf{y},\mathbf{z})\simeq C(\mathbf{x}',\mathbf{y}',\mathbf{z}')$
if any of the following conditions is satisfied:
\begin{enumerate}[label=\normalfont(\roman*)]
\setcounter{enumi}{4}
\item\label{W_5 simply_ x->x}
$\mathbf{y}'=\mathbf{y},\mathbf{z}'=\mathbf{z}$, and $\mathbf{x}'=(x_1,\ldots,x_{l-1},x_k+x_l,x_{l+1},\ldots,x_{a+n_3+n_4})$ for some distinct $k,l\in\{1,\ldots,a+n_3+n_4\}$ and $r(k)\geq r(l)$;

\item\label{W_5 simply_ y->y}
$\mathbf{x}'=\mathbf{x},\mathbf{z}'=\mathbf{z}$, and $\mathbf{y}'=(y_1,\ldots,y_{l-1},y_k+y_l,y_{l+1},\ldots,y_{a})$ for some distinct $k,l\in\{1,\ldots,a\}$ and $r(k)\leq r(l)$;

\item\label{W_5 simply_ z->z}
$\mathbf{x}'=\mathbf{x},\mathbf{y}'=\mathbf{y}$, and $\mathbf{z}'=(z_1,\ldots,z_{l-1},z_k+z_l,z_{l+1},\ldots,z_{c})$ for some distinct $k,l\in\{1,\ldots,c\}$ and $s(k)\leq s(l)$;

\item\label{W_5 simply_ y->z}
$\mathbf{x}'=\mathbf{x},\mathbf{y}'=\mathbf{y}$, and $\mathbf{z}'=(z_1,\ldots,z_{l-1},y_k+z_l,z_{l+1},\ldots,z_{c})$ for some distinct $k,l\in\{1,\ldots,c\}$ and $r(k)\leq s(l)$.
\end{enumerate}
\end{lemma}

\begin{proof}
Homotopy equivalences~\ref{eta_i has order 2} --~\ref{epsilon_k has order 2} follow immediately from~\eqref{eqn_eta xi epsilon restrictions}. (Note that in case $y_j=3$ with $r(j)=1$, then we have $3\tilde{\xi}_j\simeq\tilde{\eta}_j+\tilde{\xi}_j$, so~\ref{y_j=3} follows.) Equivalences~\ref{W_5 simply_ x->x} --~\ref{W_5 simply_ y->z} follow from Lemma~\ref{lemma_simplify mapping cone} by taking $g\colon A_k\to A_l$ to be
\begin{itemize}
\item
the identity map $A_k\to A_l$ if $A_k=A_l$
\item
the Hopf map $\eta\colon S^4\to S^3$ if $A_k=S^4$ and $A_l=S^3$
\item
the inclusion $\imath_{r_l}\colon S^3\to P^4(2^{r_l})$ if $A_k=S^3$ and $A_l=P^4(2^{r_l})$
\item
$\imath_{r_kr_l}\colon P^4(2^{r_k})\to P^4(2^{r_l})$ if $A_k=P^4(2^{r_k})$ and $A_l=P^4(2^{r_l})$
\item
$\imath_{r_l}\circ\eta$ if $A_k=S^4$ and $A_l=P^4(2^{r_l})$
\item
$\jmath_{s_ks_l}$ if $A_k=\Sigma\C\PP^4(2^{s_k})$ and $A_l=\Sigma\C\PP^2(2^{s_l})$
\item
$\jmath_{s_l}\circ\imath_{r_ks_l}$ if $A_k=P^4(2^{r_k})$ and $A_l=\Sigma\C\PP^2(2^{s_l})$
\end{itemize}
where
$\imath_{rs}\colon P^4(2^r)\to P^4(2^s)$ is the map given  in~\eqref{diagram_imath ell>k} and~\eqref{diagram_imath ell<k},
$\jmath_s\colon P^4(2^{s})\to\Sigma\C\PP^2(2^s)$ is the inclusion, and $\jmath_{rs}\colon\Sigma\C\PP^2(2^r)\to\Sigma\C\PP^2(2^s)$ is the map given in Lemma~\ref{lemma_pi_5 Sigma CP^2(2^r)}.
\end{proof}

\begin{lemma}\label{lemma_homotopy type of C_varphi for a=5}
Given the coefficients $\mathbf{x},\mathbf{y},\mathbf{z}$ of zeros and ones from~\eqref{eqn_varphi_5 simplify}, there exist
\[
u\in\{1,\ldots,a+n_3+n_4\},\qquad
v\in\{1,\ldots,a\},\qquad
w\in\{1,\ldots,c\}
\]
such that
$W_5$ is homotopy equivalent to the mapping cone of
\[
x_u\tilde{\eta}_u+y_v\tilde{\xi}_v+z_w\tilde{\epsilon}_w\colon S^5\longrightarrow W_4.
\]
Furthermore, if $y_v=1$ and $z_w=1$ then $s(w)<r(v)$.
\end{lemma}

\begin{proof}
By definition, $W_5=C(\mathbf{x},\mathbf{y},\mathbf{z})$. To prove the first part of the lemma, it suffices to show $C(\mathbf{x},\mathbf{y},\mathbf{z})\simeq C(\mathbf{x}',\mathbf{y}',\mathbf{z}')$ where
\begin{enumerate}[label=(\alph*)]
\item
$\mathbf{x}'=(0,\ldots,0,x_{u},0,\ldots,0)$
\item
$\mathbf{y}'=(0,\ldots,0,y_{v},0,\ldots,0)$
\item
$\mathbf{z}'=(0,\ldots,0,z_{w},0,\ldots,0)$.
\end{enumerate}
For~(a), let $I=\{i\mid x_i=1\}$ be the index set of those nonzero $x_i=1$ in $\mathbf{x}$. If $I$ has at most one element then we are done. Otherwise take $u$ to be the largest index in $I$ and use Lemma~\ref{lemma_htpy equiv of mapping cone of varphi_5}~\ref{eta_i has order 2} and~\ref{W_5 simply_ x->x} to annihilate other nonzero $x_i$ by $x_u$ as in the proof of Lemma~\ref{lemma_htpy equiv of mapping cone of varphi_5}.

For~(b), let $J=\{j\mid y_j=1\}$ be the index set of those nonzero $y_j$ in $\mathbf{y}$. 
If $J$ has more than one element, then take $v$ to be the smallest index in $J$ and use Lemma~\ref{lemma_htpy equiv of mapping cone of varphi_5}~\ref{xi_j has order 4 if r_j=1}, \ref{xi_j has order 2 if r_j>1} and~\ref{W_5 simply_ y->y} to annihilate other nonzero $y_j$ by $y_v$. Part~$(c)$ can be proved similarly using Lemma~\ref{lemma_htpy equiv of mapping cone of varphi_5}~\ref{epsilon_k has order 2} and~\ref{W_5 simply_ z->z}.

Next we prove the second part of the lemma. Suppose $y_v=1$ and $z_w=1$ with $s(w)\geq r(v)$. Use Lemma~\ref{lemma_htpy equiv of mapping cone of varphi_5}~\ref{epsilon_k has order 2} and~\ref{W_5 simply_ y->z} to annihilate $z_w$ by $y_v$. 
\end{proof}

We are now ready to prove the main theorem.

\begin{proof}[Proof of Theorem~\ref{thm_main theorem}]
By homotopy equivalence~\eqref{equation_Sigma X'}, we have $\Sigma X\simeq(S^2)^{\vee m}\vee W_5$.
Hence it suffices to show that $W_5$ is homotopy equivalent to $\Sigma C_f\vee\bigvee_i P_i$ where $f$ is one of the maps~\ref{type of f_S^3 v Moore}, \ref{type of f_S^2 v Moore}, \ref{type of f_two Moore}, or~\ref{type of f_one Moore}, and each $P_i$ is either a sphere, a Moore space, $\Sigma\C\PP^2$, or $\Sigma\C\PP^2(2^r)$.

By Lemma~\ref{lemma_homotopy type of C_varphi for a=5}, $W_5$ is homotopy equivalent to the mapping cone of 
\begin{equation}\label{eqn_simplest form of attaching map}
x_u\tilde{\eta}_u+y_v\tilde{\xi}_v+z_w\tilde{\epsilon}_w
\end{equation}
for some $u\in\{1,\ldots,a+n_3+n_4\},v\in\{1,\ldots,a\}$ and $w\in\{1,\ldots,c\}$, and $x_u,y_v,z_w\in\{0,1\}$ such that $s(w)<r(v)$ if $y_v=z_w=1$.
Suppose $u\leq a$. By definition $\tilde{\eta}_u,\tilde{\xi}_v$ and $\tilde{\epsilon}_w$ factor through the summands $P^4(2^{r_u}),P^4(2^{r_v})$ and $\Sigma\C\PP^2(2^{s_w})$ in $W_4$. If $u\neq v$ then the map~\eqref{eqn_simplest form of attaching map} factors through
\[
(x_u\imath_{r_u}\circ\eta^2)\bot(y_v\xi_{r_v})\bot(z_w\epsilon_{s_w})\colon S^5\to P^4(2^{r_u})^{\vee x_u}\vee P^4(2^{r_v})^{\vee y_v}\vee\Sigma\C\PP^2(2^{s_w})^{\vee z_w}.
\]
Note that this map is the suspension of $(x_u{\eta}_{r_u})\bot(y_v{\xi}_{r_v})\bot(z_w{\epsilon}_{s_w})$, so $W_5\simeq\Sigma C_f\vee\bigvee_iP_i$ with~$f$ of the form~\ref{type of f_two Moore}.
If $u=v$ then the map~\eqref{eqn_simplest form of attaching map} factors through
\[
(y_v\imath_{r_v}\circ\eta^2+y_v\xi_{r_v})\bot(z_w\epsilon_{s_w})\colon S^5\to P^4(2^{r_v})^{\vee y_v}\vee\Sigma\C\PP^2(2^{s_w})^{\vee z_w}.
\]
It is the suspension of $(y_v{\eta}_{r_v}+y_v{\xi}_{r_v})\bot(z_w{\epsilon}_{s_w})$, so $W_5\simeq\Sigma C_f\vee\bigvee_iP_i$ with $f$ of the form~\ref{type of f_one Moore}.

Suppose $u>a$. If $r(u)=\omega$ then $\tilde{\eta}_u$ factors through $S^3_u\subset W_4$ and the map~\eqref{eqn_simplest form of attaching map} factors through
\[
(x_u\eta^2)\bot(y_v\xi_{r_v})\bot(z_w\epsilon_{s_w})\colon S^5\to (S^3_u)^{\vee x_u}\vee P^4(2^{r_v})^{\vee y_v}\vee\Sigma\C\PP^2(2^{s_w})^{\vee z_w}.
\]
It is the suspension of $(x_u\eta^2)\bot(y_v{\xi}_{r_v})\bot(z_w{\epsilon}_{s_w})$, so $W_5\simeq\Sigma C_f\vee\bigvee_iP_i$ with $f$ of the form~\ref{type of f_S^2 v Moore}.
If $r(u)=\omega+1$ then $\tilde{\eta}_u$ factors through $S^4_u\subset W_4$ and the map~\eqref{eqn_simplest form of attaching map} factors through
\[
(x_u\eta)\bot(y_v\xi_{r_v})\bot(z_w\epsilon_{s_w})\colon S^5\to (S^4_u)^{\vee x_u}\vee P^4(2^{r_v})^{\vee y_v}\vee\Sigma\C\PP^2(2^{s_w})^{\vee z_w}.
\]
It is the suspension of $(x_u\eta)\bot(y_v{\xi}_{r_v})\bot(z_w{\epsilon}_{s_w})$, so $W_5\simeq\Sigma C_f\vee\bigvee_iP_i$ with $f$ of the form~\ref{type of f_S^3 v Moore}. This completes the proof.
\end{proof}



\section{Analysis of the splitting}\label{sect5}

In this section we discuss algebraic constraints on the form of the splitting given in Theorem~\ref{thm_main theorem} which arise from Poincaré duality. The case in which $X$ is a closed 5-manifold is especially tractable, since information about the decomposition of $\Sigma X$ can be read off from the easily available algebraic and geometric information contained in its characteristic classes. However, most of the analysis we provide depends only on Poincaré duality, so has wider application. We assume throughout that $X$ is an orientable Poincaré duality 5-complex, but shall place assumptions on its homology only when necessary. 

For such an $X$, the cup product pairing $H^k(X;\mathbb{Z}/2)\otimes H^{5-k}(X;\mathbb{Z}/2)\rightarrow H^5(X;\mathbb{Z}/2)$ is nondegenerate and gives rise to Wu classes $\upsilon_i(X)\in H^i(X;\mathbb{Z}/2)$ which are characterised by the equation
\[
\langle Sq(u),[X]\rangle=\langle \upsilon(X)\cup u,[X]\rangle
\]
where $u\in H^*(X;\mathbb{Z}/2)$, $Sq$ is the total Steenrod square, and $\upsilon(X)=1+\upsilon_1(X)+\dots$ is the total Wu class. Stiefel--Whitney classes $w_i(X)\in H^i(X;\mathbb{Z}/2)$ are then defined by writing
\begin{equation}\label{wuclassdefinition}
Sq(\upsilon(X))=1+w_1(X)+w_2(X)+\dots
\end{equation}

Orientability implies that $\upsilon_1(X)=w_1(X)=0$, and hence that $w_2(X)=\upsilon_2(X)$. Because $Sq^i$ vanishes in $H^*(X;\mathbb{Z}/2)$ for $i\geq3$, each of $\upsilon_3(X),\upsilon_4(X)$ and $\upsilon_5(X)$ vanishes. Consequently, \eqref{wuclassdefinition} becomes
\[
\upsilon_2(X)+Sq^1\upsilon_2(X)+Sq^2\upsilon_2(X)=w_2(X)+w_3(X)+w_4(X)+w_5(X),
\]
meaning that 
\[
w_3(X)=Sq^1w_2(X),\qquad w_4(X)=w_2(X)^2,\qquad\text{and}\qquad w_5(X)=0.
\]
Thus the entire information of the Stiefel--Whitney classes of $X$ is contained in $w_2(X)$.

When $X$ is a closed manifold, the Stiefel--Whitney classes defined above agree with the Stiefel--Whitney classes of its tangent bundle. Its Euler class vanishes, so the only additional piece of information in this case comes in the form of its first Pontryagin class $p_1(X)\in H^4(X;\mathbb{Z})$. This satisfies $\rho_2(p_1(X))=w_2(X)^2$, where $\rho_2$ is reduction mod $2$.
\begin{prop}\label{propactsteen}
Let $X$ be a Poincaré duality 5-complex. The following statements hold.
\begin{enumerate}
\item The operation
\[
Sq^2\colon H^3(X;\mathbb{Z}/2)\rightarrow H^5(X;\mathbb{Z}/2)
\]
is given by cupping with $w_2=w_2(X)$.
\item The operation 
\[
Sq^2Sq^1\colon H^2(X;\mathbb{Z}/2)\rightarrow H^5(X;\mathbb{Z}/2)
\]
is given by cupping with $w_3=w_3(X)$.
\item If $H_1(X;\Z)$ contains no 2-torsion, then the operation 
\[
Sq^2\colon H^2(X;\mathbb{Z}/2)\rightarrow H^4(X;\mathbb{Z}/2)
\]
is given by cupping with $w_2=w_2(X)$.
\end{enumerate}
\end{prop}
\begin{proof}
$(1)$ This is implied by the agreement of the second Stiefel--Whitney and Wu classes of~$X$. See equation \eqref{wuclassdefinition}.

$(2)$ For $x\in H^2(X;\mathbb{Z}/2)$ use $(1)$ to get $Sq^1(w_2\cdot x)=Sq^1Sq^2x=Sq^3x=0$. On the other hand, this gives
\[
0=Sq^1(w_2\cdot x)=w_3\cdot x+w_2\cdot Sq^1x=w_3\cdot x+Sq^2Sq^1x.
\]
We conclude that $Sq^2Sq^1x=w_3\cdot x$.

$(3)$ The cup product pairing $H^1(X;\mathbb{Z}/2)\otimes H^4(X;\mathbb{Z}/2)\rightarrow \mathbb{Z}/2$ is non-degenerate. Thus for $x\in H^2(X;\mathbb{Z}/2)$ we have that $Sq^2x$ is equal to $w_2\cdot x$ if and only if $a\cdot (Sq^2x+w_2\cdot x)=0$ for all $a\in H^1(X;\mathbb{Z}/2)$. For such an $a$ we have $Sq^2a=0$ for dimension reasons, and $Sq^1a=0$ holds because of the assumption. Thus $Sq^2(a\cdot x)=a\cdot Sq^2(x)$, and hence 
\[
a\cdot (Sq^2x+w_2\cdot x)=Sq^2(a\cdot x)+w_2\cdot (a\cdot x)=(Sq^2+w_2)\cdot a\cdot x=0,
\]
since $(Sq^2+w_2)\cdot H^3(X;\mathbb{Z}/2)=0$ by $(1)$. 
\end{proof}

As we explain below, the operations
\begin{equation}\label{topsttenrodoperations}
Sq^2\colon H^3(X;\mathbb{Z}/2)\rightarrow H^5(X;\mathbb{Z}/2),\qquad Sq^2Sq^1\colon H^2(X;\mathbb{Z}/2)\rightarrow H^5(X;\mathbb{Z}/2)
\end{equation}
give direct information about the components of the map $f$ which appears in Theorem~\ref{thm_main theorem}. On the other hand, Proposition~\ref{propactsteen} shows that the vanishing of these operations is equivalent to the vanishing of~$w_2(X),w_3(X)$, respectively. In particular, the information contained in \eqref{topsttenrodoperations} is equivalent to that contained in these Stiefel--Whitney classes, which already hold significant geometric and algebraic importance.
\begin{prop}\label{propspin5complex}
Let $X$ be an Poincaré duality 5-complex with torsion-free $H_1(X;\Z)$. In the notation of Theorem~\ref{thm_main theorem}, the map $f$ contains a component $\eta,{\xi}_r$ or ${\epsilon}_r$ for some $r$ if and only if $w_2(X)\neq0$.
\end{prop}
\begin{proof}
The point is that each of these maps is detected by $Sq^2:H^3(X;\mathbb{Z}/2)\rightarrow H^5(X;\mathbb{Z}/2)$. This is obvious for $\eta$, and it follows for ${\xi}_r$ and ${\epsilon}_r$ from the properties listed in Lemma~\ref{lemma_hmtpy gps Moore space} and Lemma~\ref{lemma_pi_5 Sigma CP^2(2^r)}. As discussed above, the $Sq^2$ is trivial if and only if $w_2(X)$ vanishes. 
\end{proof}

\begin{cor}
Let $X$ be a Poincaré duality 5-complex with torsion-free $H_1(X;\Z)$.
\begin{enumerate}
\item If $\Sigma X$ splits off a copy of $\Sigma\mathbb{CP}^2$ or $\Sigma\mathbb{CP}^2(2^r)$, then $f$ has an $\eta,{\xi}_r$ or ${\epsilon}_r$ component.
\item If $w_2(X)=0$, then no $\Sigma\mathbb{CP}^2$ or $\Sigma\mathbb{CP}^2(2^r)$ splits off $\Sigma X$.
\end{enumerate}
\end{cor}
\begin{proof}
$(1)$ A copy of $\Sigma\mathbb{CP}^2$ or $\Sigma\mathbb{CP}^2(2^r)$ splitting off $\Sigma X$ would indicate that $$Sq^2:H^2(X;\mathbb{Z}/2)\rightarrow H^4(X;\mathbb{Z}/2)$$ is nontrivial. Because of Proposition~\ref{propactsteen}, this implies that $w_2(X)\neq0$.

$(2)$ By Proposition~\ref{propactsteen}, $Sq^2:H^2(X;\mathbb{Z}/2)\rightarrow H^4(X;\mathbb{Z}/2),$
and hence $Sq^2:H^3(\Sigma X;\mathbb{Z}/2)\rightarrow H^5(\Sigma X;\mathbb{Z}/2)$, are trivial.
\end{proof}

Since $w_3(X)=Sq^1w_2(X)$, the vanishing of $w_2(X)$ implies that of $w_3(X)$. However, when~$w_2(X)\neq0$, the vanishing of $w_3(X)$ has certain implications for the topology of $X$.
\begin{prop}\label{w3nonzeroimpl}
Let $X$ be a Poincaré duality $5$-complex with $w_2(X)\neq0$. Then $f$ can contain exactly one of $\xi_1,\epsilon_1$, and this occurs if and only if $w_3(X)\neq0$. If this holds, then there is $u\in H^2(X;\mathbb{Z}/2)$ with $w_2(X)\cdot Sq^1u\neq0$, and $f$ contains $\epsilon_1$ if and only if $w_2(X)\cdot u\neq0$. If $f$ contains $\xi_1$, then it contains no $\epsilon_s$ for any $s\geq1$.
\end{prop}
\begin{proof}
Evidently, the presence of $\xi_1$ or $\epsilon_1$ appearing in $f$ is equivalent to $Sq^2Sq^1$ acting nontrivially in $H^*(X;\mathbb{Z}/2)$. By Proposition~\ref{propactsteen}, this is equivalent to $w_3(X)\neq0$. In this case there is $u\in H^2(X;\mathbb{Z}/2)$ satisfying
\[
0\neq w_3(X)\cdot u=Sq^2Sq^1u=w_2(X)\cdot Sq^1u.
\]
Now, both $\xi_1$ and $\epsilon_1$ cannot appear together, since Theorem~\ref{thm_main theorem} states that if $f$ contains both $\xi_r$ and $\epsilon_s$, then $s<r$. In particular, if $\xi_1$ appears in $f$ then no $\epsilon_s$ appears.

To distinguish between the two cases, observe that $u$ must correspond to an element in $H^*(C_f;\mathbb{Z}/2)$ which is represented by either the bottom cell of $P^3(2)$, or by the the bottom cell of $\mathbb{CP}^2(2)$. In the first case, $Sq^2u=u^2=0$, while $Sq^2u\neq0$ holds in the second. Appealing again to Proposition~\ref{propactsteen} we complete the proof of the statement.
\end{proof}

There is more that can be said in case $w_2(X)\neq0$ and $w_3(X)=0$. Let $\{(E_r^\ast(X),d_r)\}_{r\in\mathbb{N}}$ be the mod $2$ cohomology Bockstein spectral sequence for $X$. Let also $\beta_r\colon H^*(X;\mathbb{Z}/2^r)\rightarrow H^{*+1}(X;\mathbb{Z}/2)$ be the Bockstein connecting map associated with the short exact sequence 
\[
0\rightarrow\mathbb{Z}/2\rightarrow\mathbb{Z}/2^{r+1}\rightarrow\mathbb{Z}/2^r\rightarrow0.
\]
Note that if $i\colon\mathbb{Z}/2\rightarrow\mathbb{Z}/2^r$ is the inclusion, then $i_\#\beta_r=d_r$ under the identification of $d_r$ given in \cite[Proposition 10.4]{McCleary:2001}. In particular, $\beta_1=Sq^1$, so $w_3(X)=0$ is the statement that $w_2(X)$ survives past the first page of the Bockstein spectral sequence.

Now assume that $w_2$ survives to the $r^{th}$ page of the Bockstein spectral sequence. Then there is a class $w^{(r)}_2\in H^2(X;\mathbb{Z}/2^r)$ with $\rho_2w^{(r)}_2=w_2$. Because $H^5(X;\mathbb{Z})$ is free abelian, $\beta_r(H^4(X;\mathbb{Z}/2^r))=0$, so for any $x\in H^2(X;\mathbb{Z}/2^r)$ it holds that 
\[
0=\beta_r(w_2^{(r)}\cdot x)=\beta_r(w_2^{(r)})\cdot \rho_2(x)+w_2\cdot\beta_rx=\beta_r(w_2^{(r)})\cdot \rho_2(x)+Sq^2\beta_r(x).
\]
That is,
\[
Sq^2\beta_r(x)=\beta_r(w_2^{(r)})\cdot \rho_2(x).
\]
Note that the right-hand side of this equation is independent of the choice of $w_2^{(r)}$. In any case, if $\beta_r(w_2^{(r)})=0$, then $Sq^2\beta_r$ acts trivially on $H^*(X;\mathbb{Z}/2^r)$. On the other hand, if $d_rw_2\neq0$, then necessarily $\beta_r(w_2^{(r)})\neq0$, and the operation $Sq^2\beta_r$ is nontrivial. From this discussion we have the following.
\begin{prop}\label{w3zeroimplication}
Let $X$ be a Poincaré duality 5-complex with $H_1(X;\Z)$ torsion-free and $w_2(X)\neq0$. Write $r(X)\in\mathbb{N}\cup\{\infty\}$ for the greatest integer for which $w_2(X)$ survives to the $E_{r(X)}$-page of the mod $2$ Bockstein spectral sequence. A necessary and sufficient condition that $r(X)>1$ is that $w_3(X)=0$. The following statements hold.
\begin{enumerate}
\item $r(X)<\infty$ if and only if exactly one of $\xi_{r(X)},\epsilon_{r(X)}$ appears in $f$, but no $\xi_s$ or $\epsilon_s$ for $s<r(X)$ does. If $\epsilon_{r(X)}$ appears, then there is $u\in H^2(X;\mathbb{Z}/2^{r(X)})$ with $Sq^2\beta_{r(X)}(u)\neq0\neq Sq^2\rho_2(u)$. If $\xi_{r(X)}$ appears, then no $\epsilon_s$ appears for any $s>r(X)$.
\item $r(X)=\infty$ if and only if $f\simeq\eta$. If this holds, then there is $u\in H^3(X;\mathbb{Z})$ such that $Sq^2\rho_2(u)\neq0$.
\end{enumerate}
\end{prop}
\begin{proof}
Only the finer points need be explained. In $(1)$ these are dealt with by appealing to the part of Theorem~\ref{thm_main theorem} that states that if $f$ contains both $\xi_r$ and $\epsilon_s$, then $s<r$. 

For $(2)$, it has been shown above that $f$ contains no $\xi_r$ or $\epsilon_s$ for any $r,s\geq1$. However, Proposition~\ref{propspin5complex} then states that $f$ must have an $\eta$ component. There are four possible options for $f$ listed in Theorem~\ref{thm_main theorem}, and given the constraints just mentioned, $f\simeq\eta$ is the only one which may be realised.
\end{proof}

Next, components $\eta^2,\iota_r\circ\eta^2$ appearing in $f$ can be detected using secondary cohomology operations. In particular, this is possible using the secondary operation based on the relation $Sq^2Sq^2+Sq^1(Sq^2Sq^1)=0$ \cite[p.149]{harper}. Assuming that $H_1(X;\Z)$ contains no 2-torsion, Proposition~\ref{propactsteen} may be used to construct this operation as
\[
\Theta\colon\{x\in H^2(X;\mathbb{Z}/2)\mid x\cdot w_2=0=x\cdot w_3\}\rightarrow H^5(X;\mathbb{Z}/2)/(w_2\cdot H^3(X;\mathbb{Z}/2)).
\]
Unfortunately, use of this operation is only feasible when $w_2=0$. In this case it is an operation
\begin{equation}\label{secoper}
\Theta\colon H^2(X;\mathbb{Z}/2)\rightarrow H^5(X;\mathbb{Z}/2).
\end{equation}
Alone it cannot distinguish between $\eta^2,\iota_r\circ\eta^2$, but this may be accomplished by means of a higher Bockstein operator $\beta_r$.
\begin{prop}
Let $X$ be a Poincaré duality 5-complex which has torsion-free $H_1(X)$ and $w_2(X)=0$. If~$f$ has a component of either $\eta^2$ or $\iota_r\circ\eta^2$, then there is $u\in H^2(X;\mathbb{Z}/2)$ satisfying $\Theta(u)\neq0$. If $\beta_ru=0$ for all $r$, then $u$ detects $\eta^2$. If $\beta_ru\neq 0$ for some $r\geq2$, then $u$ detects $\iota_r\circ\eta^2$. \qed
\end{prop}
Combining the operations above makes it possible also to detect a $\iota_r\circ\eta^2+\xi_r$ ($r\geq2$) component of $f$.


\begin{proof}[Proof of Theorem~\ref{corollary1}]
$(1)$ We refer to Theorem~\ref{thm_main theorem}. Proposition~\ref{propspin5complex} explains that $f$ can only contain components $\eta^2$ or $\imath_r\circ\eta^2$, both of which are detected by the secondary operation~$\Theta$ of equation~\eqref{secoper}. Consequently, the complex $\Sigma C_f$ has at most three cells, and splits up further if and only if $\Theta$ evaluates trivially on its cohomology. With this reduction, it is easy to follow the method of \cite[p.32]{MadsenMilgram:1979} to see that $X$ being a manifold forces $\Theta$ to be trivial.

$(2)$ The first part of the statement is covered by Proposition~\ref{w3nonzeroimpl}. For the second, observe that $w_2(X)\cdot w_3(X)=w_2(X)\cdot Sq^1w_2(X)\neq0$. Thus, still following Proposition~\ref{w3nonzeroimpl}, $\xi_1$ is detected if $w_2(X)^2=0$. However, the relation $w_2(X)^2=\rho_2(p_1(X))$ holds for any 5-dimensional manifold. Clearly $\rho_2(p_1(X))=0$ if and only if $p_1(X)$ is divisible by $2$.

$(3)$ See Proposition~\ref{w3zeroimplication}.
\end{proof}
\begin{eg}
Let $S_g$ be a closed orientable surface of genus $g$, and $Y$ a closed, orientable~3-manifold with $H_1(Y;\Z)$ torsion-free. Form $X=S_g\times Y$. Then $X$ is spin with torsion-free homology. Theorem~\ref{corollary1} says that
\[
\Sigma X\simeq \Sigma S_g\vee \Sigma Y\vee \Sigma(S_g\wedge Y)
\]
splits as a wedge of spheres. Clearly this requires that $\Sigma Y$ is a wedge of spheres. \qed
\end{eg}
\begin{eg}
Let $Y$ be a nonspin 4-manifold with $H_1(Y;\Z)\cong\mathbb{Z}^{m-1}$ and let $X=S^1\times Y$. Then
\[
\Sigma X\simeq S^2\vee \Sigma Y\vee \Sigma^2Y.
\]
The homology of $X$ is torsion free, so part $(3)$ of Theorem~\ref{corollary1} is in effect with $w_2(X)$ surving to the $E_\infty$-page of the Bockstein spectral sequence. Thus $f\simeq\eta$. This should be compared with the decomposition of $\Sigma X$ given in \cite[Theorem 1.1]{ST19}, which implies the same result. \qed
\end{eg}


For simply connected $X$, the maps $f$ have been determined completely by St\"ocker \cite{Stocker:1982}. In the case that $X$ is a smooth 5-manifold, much of this was already implicit in Barden's classification results \cite{Barden:1965}. However, St\"ocker's results also extend to simply connected Poincar\'e duality complexes. We end this section by reviewing these known results.

In \cite{Barden:1965} Barden defines simply connected, oriented 5-manifolds $X_0,X_{-1},X_{\infty},M_\infty$, and $X_k,M_k$ for $k\geq2$, and proves the following.

\begin{thm}[Barden \cite{Barden:1965}, Theorem 2.3]\label{Bardenclass}
Every simply connected, closed, orientable, smooth 5-manifold is diffeomorphic to a manifold of the form
\[
X_{j,k_1,\dots,k_n}=X_j\# M_{k_1}\#\cdots\# M_{k_n}
\]
where $-1\leq j\leq\infty$ and $1<k_1\leq k_2\leq\dots\leq k_n$ are such that either $k_i$ divides $k_{i+1}$ or~$k_{i+1}=\infty$. \qed
\end{thm}

We have $X_0=S^5$ and $M_\infty=S^2\times S^3$. The manifolds $M_k$ are of the form $$M_k\simeq(P^2(k)\vee P^2(k))\cup_\omega e^5,$$
where $\omega$ is a Whitehead product \cite{Stocker:1982}. All of these manifolds split apart after a single suspension. The other manifolds are more interesting.
\begin{itemize}
\item $X_{-1}=SU_3/SO_3$ is the Wu manifold. Jie Wu \cite[Example 6.15]{Wu:2003} has observed that
\[
X_{-1}\simeq P^3(2)\cup_{\xi_1} e^5.
\]
\item $X_\infty=S^2\widetilde\times S^3$ is the total space of the nontrivial $S^3$-bundle over $S^2$, classified by a generator of $\pi_1(SO_4)\cong\mathbb{Z}/2$. Following \cite{JamesWhitehead:1954} we check that
\[
X_\infty\simeq (S^2\vee S^3)\cup_{\phi_{\infty}} e^5,
\]
where $\phi_{\infty}=\jmath_3\circ\eta+[\iota_2,\iota_3]$. 
\item $X_k$ for $k\geq1$ is a nonspin manifold which, according to \cite{Stocker:1982}, is of the form
\[
X_k\simeq(P^3(2^k)\vee P^3(2^k))\cup_{\phi_k} e^5,
\]
where $\phi_k=\omega+\jmath_1\circ\xi_{k}$ with $\omega$ a Whitehead product.
\end{itemize}

In addition to these manifolds there are certain interesting non-smoothable, simply connected Poincar\'e duality 5-complexes. One such is the Gitler-Stasheff example \cite{GitlerStasheff:1965}
\[
M'_\infty\simeq (S^2\vee S^3)\cup_\varphi e^5,
\]
where $\varphi=[\iota_2,\iota_3]+\jmath_2\circ\eta^2$. The others we consider are those constructed by St\"ocker \cite{Stocker:1982};
\begin{itemize}
\item $M'_k\simeq(P^2(k)\vee P^2(k))\cup_{\varphi'_k} e^5$ for $k\geq2$, where $\varphi'_k=\omega+\iota_k\circ\eta^2$ with $\omega$ a Whitehead product.
\item $X_k'=(P^3(2^k)\vee P^3(2^k))\cup_{\phi'_k}e^5$ for $k\geq1$, where $\phi'_k=\omega+\jmath_1(\xi_{k}+\iota_{2^k}\circ\eta^2)$ with $\omega$ a Whitehead product.
\end{itemize}
\begin{thm}[St\"ocker \cite{Stocker:1982}, Theorem 10.1]
Every simply connected, 5-dimensional, orientable Poincar\'e duality complex is homotopy equivalent to a connected sum of the spaces defined above. \qed
\end{thm}

\begin{cor}\label{corollarysimplyconnectclass}
If $X$ is a closed, orientable, simply connected 5-manifold, then $f$ can contain only $\eta$ or $\xi_r$ for some $r\geq1$. If $X$ is an orientable, simply connected Poincar\'e duality 5-complex, then $f$ can contain only $\eta,\eta^2,{\xi}_r,\iota_r\circ\eta^2$, or $({\xi}_r+\iota_r\circ\eta^2)$ for some $r\geq1$. \qed
\end{cor}

\section{Applications}\label{sect6}

\subsection{Generalised cohomology}\label{sectgencohom}

As a simple application of our main theorem we explain how to evaluate a given cohomology theory $h^*$ on any Poincar\'{e} duality complex satisfying the assumptions of Theorem~\ref{thm_main theorem}. For example $h^*$ could be singular cohomology, complex or real $K$-theory, or cobordism. We will give applications in the sequel when $h^*$ is stable cohomotopy.

\begin{thm}\label{Thm51}
Let $X$ be a connected, orientable, 5-dimensional Poincar\'{e} duality complex with homology as in~(\ref{table_original M hmlgy}), and let $h^*$ be a reduced generalised cohomology theory. Then 
\begin{align*}
h^*(X)\cong&\;\textstyle{\bigoplus^mh^*(S^1)\oplus\bigoplus^{n_3} h^*(S^2)\oplus\bigoplus^{n_4}h^*(S^3)\oplus\bigoplus^{n_5} h^*(S^4)\oplus}\\
&\;\textstyle{\bigoplus_{t'_i\in T'} h^*(P^3(t'_i))\oplus\bigoplus^bh^*(\mathbb{CP}^2)\oplus\bigoplus^c_{j=1}h^*(\mathbb{CP}^2(2^{r_j}))\oplus h^*(C_f)}
\end{align*}
where the integers $n_3,n_4,n_5,a,b,c,t'_i,r_j$ are explained in Theorem~\ref{thm_main theorem}.
\end{thm}

\begin{proof}
The lemma follows immediately from the string of isomorphisms
\[
\textstyle h^*(X)\cong h^{\ast+1}(\Sigma X)\cong h^{*+1}(\bigvee_i\Sigma P_i)\cong\bigoplus_ih^{*+1}(\Sigma P_i)\cong\bigoplus_i h^{*}(P_i).
\]
where the first and the third isomorphisms are due to the suspension isomorphisms of $h^*$ and the second isomorphism is due to the suspension splitting $\Sigma X\simeq\bigvee_i\Sigma P_i$ in Theorem~\ref{thm_main theorem}.
Since we work only with finite wedge sums, we do not need to assume that $h^*$ satisfies the wedge axiom to obtain the third isomorphism.
\end{proof}

There are similar decompositions of $h_*(M)$ for any generalised homology theory $h_*$ which we will leave to the reader to spell out. Interesting applications of this are when $h_*$ is~$K$-homology or stable homotopy.

\subsection{Cohomotopy groups and sets}\label{sectcohomgroups}

In this section we consider the problem of computing the cohomotopy sets $\pi^n(X)=[X,S^n]$ and $\pi^n(X;\mathbb{Z}/k)=[X,P^{n+1}(k)]$ for $X$ a Poincar\'{e} duality 5-complex satisfying the assumptions of our main theorem. We start by briefly recalling the definitions of these sets and the construction of group operations on them when~$X$ is any CW complex.

In general, $\pi^n(X)$ is only a pointed set. However, if $X$ is a CW complex of dimension~\mbox{$\leq 2n-2$}, then $\pi^n(X)$ carries a natural abelian group structure, which can defined as follows. Suppose maps $f,g\colon X\rightarrow S^n$ are given. Because the inclusion $S^n\vee S^n\hookrightarrow S^n\times S^n$ is $(2n-1)$-connected, there is a unique compression of $(f,g)\colon X\rightarrow S^n\times S^n$ into a map $\theta_{f,g}\colon X\rightarrow S^n\vee S^n$. We define
\begin{equation}\label{dfn_cohtmpy gp structure 1}
f+g=\nabla\circ\theta_{f,g}\colon X\rightarrow S^n,
\end{equation}
where $\nabla\colon S^n\vee S^n\rightarrow S^n$ is the folding map, and put $-f=(-1)\circ f$. Borsuk \cite{Borsuk:1936} shows that these definitions equip $\pi^n(X)$ with an abelian group structure (see Spanier \cite{Spanier:1949} for full proofs).

There is also a second description of this group structure which is due to Taylor \cite{Taylor:2012}, and can be useful to know. If $\dim X\leq2n-2$, then by a standard argument we have $\pi^n(X)\cong[X,PK_{2n-2}(S^n)]$, where $PK_i(S^n)$ denotes the $i$th Postnikov section of $S^n$. Owing essentially to the Freudenthal Suspension Theorem, $PK_{2n-2}(S^n)$ is an infinite loop space. Hence $\pi^n(X)$ is an abelian group in each of these cases. 

This viewpoint has some advantages. For example, it makes clear that while working in this dimension range, the function 
\begin{equation}\label{cohomsusphom}
\Sigma\colon \pi^n(X)\rightarrow\pi^{n+1}(\Sigma X)
\end{equation}
induced by the suspension map $\sigma\colon S^n\rightarrow\Omega S^{n+1}$ is compatible with group operations on both sides. This gives a well-defined suspension {\it homomorphism}.

In fact, this leads to yet another description of the group $\pi^n(X)$. Namely, the suspension~\eqref{cohomsusphom} is bijective when $\dim X\leq 2n-2$, and we may turn $\pi^n(X)$ into a group by requiring it to be a homomorphism. We've already explained why this gives the same operation as Taylor's. Showing that it agrees with Borsuk's original definition~\eqref{dfn_cohtmpy gp structure 1} is the easiest way to establish the equivalence of all three definitions.

Now, Borsuk's original definition~\eqref{dfn_cohtmpy gp structure 1} has its merits. In case $n=1,3,7$ there is one further group structure on $\pi^n(X)$, which comes from the H-space multiplication on $S^n$. Because the multiplication represents an extension of the folding map $\nabla$ over $S^n\times S^n$, we find further agreement between this operation and the three set out above.

It is also possible to define cohomotopy groups with coefficients. These were introduced originally by Peterson \cite{Peterson:1956} via the definition
\[
\pi^n(X;\mathbb{Z}/k)=[X,P^{n+1}(k)].
\]
As above, $\pi^n(X;\mathbb{Z}/k)$ is a group when $X$ is a CW complex of dimension $\leq 2n-2$. Moreover, the maps $\pi^n(X)\to\pi^n(X;\mathbb{Z}/k)$ and $\pi^n(X;\mathbb{Z}/k)\rightarrow \pi^{n+1}(X)$ which are induced by the inclusion $S^n\hookrightarrow P^{n+1}(k)$, and the pinch map $P^{n+1}(k)\to S^{n+1}$, respectively, are group homomorphisms in this range. More recently, these groups have been studied by Li, Pan, and Wu \cite{LiPanWu:2023}, who called them {\it modular cohomotopy groups}.

The cohomotopy and modular cohomotopy groups are related in the usual long exact sequence which starts 
\[
\pi^m(X)\xrightarrow{\times k}\pi^m(X)\rightarrow\pi^m(X;\mathbb{Z}/k)\rightarrow\pi^{m+1}(X)\xrightarrow{\times k}\pi^{m+1}(X)\rightarrow\dots
\]
and extends infinitely to the right, where $m\geq\max\{(\dim X+1)/2,3\}$ (these sets may fail to be groups in the bottom degree). For example, when $\dim X=5$, the following sequence is exact
\[
\pi^3(X)\xrightarrow{\times k}\pi^3(X)\rightarrow\pi^3(X;\mathbb{Z}/k)\rightarrow\pi^4(X)\rightarrow\dots\rightarrow\pi^5(X)\xrightarrow{\times k}\pi^5(X)\rightarrow\pi^5(X;\mathbb{Z}/k)\rightarrow0.
\]

The following statement is a simple application of the Freudenthal Theorem which will be used in the sequel.
\begin{prop}\label{cohomisstablecohom}
If $X$ is a CW complex of dimension $\leq 2n-2$, then stabilisation induces isomorphisms
\[
\pi^n(X)\cong\pi^n_S(X),\qquad\qquad \pi^n(X;\mathbb{Z}/k)\cong \pi^n_S(X;\mathbb{Z}/k)
\]
where $\pi^n_S(X)$ ($\pi^{n+1}_S(X;\mathbb{Z}/k)$) is the nth stable cohomotopy group (with coefficients) of $X$. \qed
\end{prop}
Now, let us turn to computation. For this we specialise to $X$ being a 5-dimensional Poincar\'{e} duality complex satisfying the assumptions of Theorem~\ref{thm_main theorem}. Of course, 
\[
\pi^n(X)=0=\pi^n(X;\mathbb{Z}/k)\qquad n\geq6
\]
for dimension reasons. We also compute
\[
\pi^5(X)\cong H^5(X;\Z)\cong \mathbb{Z},\qquad \pi^5(X;\mathbb{Z}/k)\cong H^5(X;\mathbb{Z}/k)\cong \mathbb{Z}/k
\]
using the Hopf theorem. In this range, $\pi^4(X)$ and $\pi^4(X;\Z/k)$ are also groups, and their structure is described below using results known in the literature.

On the other hand, $\pi^3(X)$ and $\pi^3(X;\Z/k)$ fall outside of Borsuk's range, and must be dealt with separately. While the group structure on $S^3$ induces a group structure on $\pi^3(X)$, a priori $\pi^3(X;\mathbb{Z}/k)$ is only a pointed set. As it turns out, $\pi^3(X;\mathbb{Z}/k)$ does carry a canonical group structure. This, along with the group $\pi^3(X)$, is described in Sections \ref{sectpi3X} and \ref{sectpi3xzk} below. We also give information about the set $\pi^2(X)$ in Section~\ref{sectpi2X}.

As for $\pi^4(X)$, we have a well-known result of Steenrod \cite{Steenrod:1947}, which gives a short exact sequence\begin{equation}\label{steenrodexactseq}
0\rightarrow H^5(X;\mathbb{Z}/2)/Sq^2(H^3(X;\mathbb{Z}))\rightarrow \pi^4(X)\rightarrow H^4(X;\Z)\rightarrow0.
\end{equation}
Under our standing assumption that $H_1(X;\Z)\cong H^4(X;\Z)$ is torsion-free, the sequence splits (as observed by Taylor \cite[Example 6.3]{Taylor:2012}, this splitting holds under slightly more general conditions). We can also replace $H^3(X;\mathbb{Z})$ with the mod $2$ reduced group. Thus
\begin{prop}\label{pi4X}
There is a group isomorphism
\[
\pi^4(X)\cong H^4(X;\Z)\oplus H^5(X;\mathbb{Z}/2)/Sq^2(H^3(X;\mathbb{Z}/2)). \qed
\]
\end{prop}

In case $X$ is a closed 5-manifold, this may be written
\[
\pi^4(X)\cong \begin{cases}
H^4(X;\Z)&\text{if $X$ is nonspin}\\
H^4(X;\Z)\oplus\mathbb{Z}/2&\text{if $X$ is spin}.
\end{cases}
\]
This situation has been studied in detail by Konstantis \cite{Konstantis:2021}, who has explicitly described a splitting map for Steenrod's exact sequence \eqref{steenrodexactseq}. Konstantis gets some mileage out of the fact that $\pi^4(X)\cong[X,\mathbb{HP}^\infty]$, and elements of the latter set are in one-to-one correspondence with isomorphism classes of quaternionic line bundles over $X$.

To introduce coefficients, consider the exact sequence
\[
\dots\rightarrow \pi^4(X)\xrightarrow{\times k}\pi^4(X)\rightarrow\pi^4(X;\mathbb{Z}/k)\rightarrow \pi^5(X)\xrightarrow{\times k}\pi^5(X)\rightarrow\dots
\]
Since multiplication by $k$ on the torsion-free group $\pi^5(X)\cong\mathbb{Z}$ is injective, the map $\pi^4(X)\rightarrow\pi^4(X;\mathbb{Z}/k)$ is onto. Hence $\pi^4(X;\mathbb{Z}/k)\cong \pi^4(X)/(k\cdot \pi^4(X))$. Making use of Proposition~\ref{pi4X} we thus state the following.

\begin{prop}
For an odd integer $k\geq1$ it holds that
\[
\pi^4(X;\mathbb{Z}/k)\cong H^4(X;\mathbb{Z}/k).
\]
For an even integer $k\geq2$ it holds that
\[
\pi^4(X;\mathbb{Z}/k)\cong H^4(X;\mathbb{Z}/k)\oplus H^5(X;\mathbb{Z}/2)/Sq^2(H^3(X;\mathbb{Z}/2)).
\]
\end{prop}

The reader may wish to compare these with the results in \cite{LiPanWu:2023}, where for a prime $p$ and integer $r\geq1$ the short exact sequence
\[
0\rightarrow H^{n+2p-3}(Y;\mathbb{Z}/p)/\mathcal{P}^1(H^{n-1}(Y;\mathbb{Z}/p^r))\rightarrow \pi^n(Y;\mathbb{Z}/p^r)\rightarrow H^n(Y;\mathbb{Z}/p^r)\rightarrow 0
\]
is constructed for any complex $Y$ with $\dim Y\leq n+2p-3$ and any $n\geq2p-1$.

\subsubsection{\textbf{The group $\pi^3(X)$}.}\label{sectpi3X} The group structure on $\pi^3(X)$ comes from the Lie multiplication on $S^3$, as it is outside of the dimension range of Borsuk's constructions. For the same reason, Proposition \ref{cohomisstablecohom} is not immediately applicable.


\begin{prop}\label{proppi3X}
There is a group isomorphism
\begin{align*}
\pi^3(X)\cong \pi^3_S(X)\cong& {\textstyle \pi_S^3(S^3)^{\oplus n_4}\oplus\pi_S^3(S^4)^{\oplus{n_5}}\oplus \bigoplus_{t'_i\in T'} \pi^3_S(P^3(t'_i))\oplus\pi_S^3(\mathbb{CP}^2)^{\oplus b}\oplus}\\
&{\textstyle \bigoplus^c_{j=1}\pi^3_S(\mathbb{CP}^2(2^{r_j}))\oplus \pi_S^3(C_f)}
\end{align*}
where
\begin{itemize}
\item $\pi^3_S(S^3)\cong\mathbb{Z}$ and $\pi^3_S(S^4)\cong\mathbb{Z}/2$,
\item $\pi^3_S(P^3(t))\cong\mathbb{Z}/t$, 
\item $\pi^3_S(\mathbb{CP}^2)=0$, 
\item $\pi^3_S(\mathbb{CP}^2(2^r))\cong\mathbb{Z}/2^{r+1}$.
\end{itemize}
\end{prop}
\begin{proof}
Since $S^3$ is an H-space, the suspension map $\sigma\colon S^3\rightarrow \Omega S^4$ has a left homotopy inverse. Consequently, the induced map
\begin{equation}\label{suspensionindim3}
\Sigma\colon\pi^3(X)\rightarrow\pi^4(\Sigma X)
\end{equation}
is injective. On the other hand, its surjectivity is already covered by the Freudenthal Theorem. Thus \eqref{suspensionindim3} is bijective.

Now, although $\sigma$ is not an H-map, the obstruction to it being one is represented by a mapping $S^3\wedge S^3\cong S^6\rightarrow \Omega S^4$. Since $X$ is 5-dimensional, the obstruction vanishes when applied to any mapping $X\rightarrow S^3$. It follows that \eqref{suspensionindim3} is an isomorphism of groups.

Finally, $\Sigma X$ is 6-dimensional, so Proposition \ref{cohomisstablecohom} may now be applied to $\pi^4(\Sigma X)$ to complete the proof.

As for the bullet points, the groups $\pi^3_S(S^n)=\pi_n^S(S^3)$ are well-known, and $\pi^3_S(P^3(t))\cong\mathbb{Z}/t$ and $\pi^3_S(\mathbb{CP}^2)=0$ are easily obtained.

For $\pi^3_S(\mathbb{CP}^2(2^{r}))$, apply $\pi^3_S$ to the cofibre sequence $S^3\rightarrow P^3(2^r)\rightarrow \mathbb{CP}^2(2^r)$. The first map factors through the inclusion $S^2\hookrightarrow P^3(2^r)$, so induces the trivial map on both $\pi^3_S$ and $\pi^4_S$. Thus we obtain a short exact sequence
\begin{equation}\label{shortxactseqforpi3cp22r}
0\leftarrow \mathbb{Z}/2^r\leftarrow \pi^3_S(\mathbb{CP}^2(2^r))\leftarrow\mathbb{Z}/2\leftarrow 0.
\end{equation}
On the other hand, there is a cofibre sequence $\mathbb{CP}^2\rightarrow\mathbb{CP}^2(2^r)\rightarrow S^3$ whose exact sequence gives an epimorphism $\pi^3_S(S^3)\cong\mathbb{Z}\rightarrow \pi^3_S(\mathbb{CP}^2(2^r))$. It follows from this that $\pi^3_S(\mathbb{CP}^2(2^r))$ is generated by a single element, and hence that \eqref{shortxactseqforpi3cp22r} cannot split. But $\text{Ext}(\mathbb{Z}/2^r,\mathbb{Z}/2)\cong\mathbb{Z}/2$, so the only other possibility is $\pi^3_S(\mathbb{CP}^2(2^r))\cong\mathbb{Z}/2^{r+1}.
$
\end{proof}
We do not give detailed information about the group $\pi^3_S(C_f)$, since it depends so much on the exact form of the attaching map $f$.

\subsubsection{\textbf{The set $\pi^3(X;\mathbb{Z}/k)$}.}\label{sectpi3xzk} Recall that the suspension map
\[
\sigma\colon P^{n+1}(k)\rightarrow \Omega P^{n+2}(k)
\]
is $(2n-1)$-connected. Thus $\Sigma\colon \pi^n(X;\mathbb{Z}/k)\rightarrow\pi^{n+1}(\Sigma X;\mathbb{Z}/k)$ is bijective for $\dim X\leq 2n-2$ and is surjective for $\dim X=2n-1$. This was partially recorded in Proposition \ref{cohomisstablecohom} as the statement that if $X$ is a 5-dimensional CW complex, then $\Sigma\colon\pi^n(X;\mathbb{Z}/k)\rightarrow\pi^{n+1}(\Sigma X;\mathbb{Z}/k)$ is an isomorphism for $n\geq4$. The following extends this range by one dimension.

\begin{prop}\label{pi3withcoeffsisagroup}
Let $X$ be a 5-dimensional CW complex. Then for any $k\geq1$ the suspension $\Sigma\colon\pi^3(X;\mathbb{Z}/k)\rightarrow\pi^4(\Sigma X;\mathbb{Z}/k)$ is a bijection.
\end{prop}

\begin{proof}
Write $P^q=P^q(k)$ and identify $\Omega\Sigma P^4$ with the James construction $J(P^4)$. Then the suspension $\sigma$ is identified with the inclusion $P^4=J_1(P^4)\hookrightarrow J(P^4)$. Letting $J_2(P^4)$ denote the second stage of the canonical filtration on $J(P^4)$, we have $[X,J(P^4)]\cong [X,J_2(P^4)]$ for dimension reasons. Furthermore, it is well-known that this space sits in a cofibration
\[
\Sigma P^3\wedge P^3\xrightarrow{[\one,\one]}P^4\rightarrow J_2(P^4)
\]
where $[\one,\one]$ is the Whitehead product $\nabla\circ[1,1]$. Thus the $6$-skeleton of $J_2(P^4)$ is of the form $P^4\cup_\alpha e^6$, where $\alpha$ is the restriction of $[\one,\one]$ to the bottom cell of $\Sigma P^3\wedge P^3$. Let $\imath\colon S^3\rightarrow P^4$ be the inclusion and consider the diagram 
\[
\xymatrix{\Sigma S^2\wedge S^2\ar[r]^-{\Sigma\imath\wedge\imath}\ar[d]_-{[\one,\one]}&\Sigma P^3\wedge P^3\ar[d]^-{[\one,\one]}\\
S^3\ar[r]^-\imath&P^4.}
\]
The diagram homotopy commutes by naturality of the Whitehead product and shows that~$\alpha$ factors through $[\one,\one]\colon\Sigma S^3\wedge S^3\rightarrow S^3$. However, since $S^3$ is an H-space, this latter Whitehead product is trivial. The $6$-skeleton of $J_2(P^4)$ is then given by $P^4\vee S^6$, and it follows that the induced map
\[
\sigma_*:[X,P^4]\rightarrow [X,J(P^4)]\cong [X,P^4\vee S^6]\cong [X,P^4]
\]
is bijective.
\end{proof}

Proposition~\ref{cohomisstablecohom} gives $\pi^4(\Sigma X;\mathbb{Z}/k)\cong \pi^4_S(\Sigma X;\mathbb{Z}/k)$. Combining this with Proposition~\ref{pi3withcoeffsisagroup}, we record the following.

\begin{cor}
If $X$ is a 5-dimensional CW complex, then for any $k\geq1$, the set
\[
\pi^3(X;\mathbb{Z}/k)\cong \pi^3_S(X;\mathbb{Z}/k)
\]
carries a canonical group structure.
\end{cor}

Finally, Theorem~\ref{Thm51} yields a decomposition result for this group.
\begin{prop}
Let $X$ be as in Theorem~\ref{thm_main theorem}. Then for any $k\geq1$ there is a group isomorphism
\begin{align*}
\pi^3(X;\mathbb{Z}/k)\cong&\;\textstyle{\pi^3_S(S^3;\mathbb{Z}/k)^{\oplus n_4}\oplus\pi^3_S(S^4;\mathbb{Z}/k)^{\oplus n_5}\oplus\bigoplus_{t'_i\in T'} \pi^3_S(P^3(t'_i);\mathbb{Z}/k)\oplus}\\
&\;\textstyle{\pi^3_S(\mathbb{CP}^2;\mathbb{Z}/k)^{\oplus b}\oplus\bigoplus^c_{j=1}\pi^3_S(\mathbb{CP}^2(2^{r_j});\mathbb{Z}/k)\oplus \pi^3_S(C_f;\mathbb{Z}/k)}
\end{align*}
where
\begin{itemize}
\item $\pi^3_S(S^3;\mathbb{Z}/k)=\mathbb{Z}/k$,
\item $\pi^3_S(S^4;\mathbb{Z}/k)=\begin{cases}0&k\;\text{odd}\\\mathbb{Z}/2&k\;\text{even,}\end{cases}$
\item $\pi^3_S(P^3(t);\mathbb{Z}/k)\cong \mathbb{Z}/gcd(t,k)$,
\item $\pi^3_S(\mathbb{CP}^2;\mathbb{Z}/k)=0$,
\item $\pi^3_S(\mathbb{CP}^2(2^r);\mathbb{Z}/k)=\mathbb{Z}/gcd(2^{r+1},k)$.
\end{itemize}
\end{prop}
\begin{proof}

We compute the component groups. The first two are given in Lemma~\ref{lemma_hmtpy gps Moore space}, and for the third, we have $\pi^3_S(P^3(t);\mathbb{Z}/k)=\pi_3^S(P^4(t);\mathbb{Z}/k)$. A universal coefficient theorem shows that $\pi_3^S(P^4(t);\mathbb{Z}/k)\cong\pi_3^S(P^4(t))\otimes\mathbb{Z}/k$, which gives the claimed result.

As for the fourth group, there is an exact sequence
\[
0\leftarrow\pi^3_S(\mathbb{CP}^2;\mathbb{Z}/k)\leftarrow \pi^3_S(S^4;\mathbb{Z}/k)\xleftarrow{\eta^*} \pi^3_S(S^3;\mathbb{Z}/k)\leftarrow\dots
\]
The group $\pi^3_S(S^4;\mathbb{Z}/k)$ vanishes if $k$ is odd. If $k$ is even, then it is $\mathbb{Z}/2$, and the generators given in Lemma~\ref{lemma_hmtpy gps Moore space} show that the $\eta^*$ appearing in the sequence is onto.

For the last group use the exact sequence
\[
\dots\rightarrow\pi^3(\mathbb{CP}^2(2^r))\xrightarrow{\times k}\pi^3(\mathbb{CP}^2(2^r))\rightarrow \pi^3_S(\mathbb{CP}^2(2^r);\mathbb{Z}/k)\rightarrow \pi^4(\mathbb{CP}^2(2^r))\xrightarrow{\times k}\pi^4(\mathbb{CP}^2(2^r))\rightarrow\dots
\]
Since $\pi^4(\mathbb{CP}^2(2^r))\cong\mathbb{Z}$, the right-most arrow is injective. On the left-hand side, it was shown in Proposition~\ref{proppi3X} that $\pi^3(\mathbb{CP}^2(2^r))\cong\mathbb{Z}/2^{r+1}$. This means that
\[
\dots\rightarrow\mathbb{Z}/2^{r+1}\xrightarrow{\times k}\mathbb{Z}/2^{r+1}\rightarrow \pi^3_S(\mathbb{CP}^2(2^r);\mathbb{Z}/k)\rightarrow0
\]
is exact, and this gives the result immediately.
\end{proof}

\subsubsection{\textbf{The set $\pi^2(X)$}.}\label{sectpi2X} There is no natural group structure on $\pi^2(X)$. However, this set does carry a $\pi^3(X)$-action which contains useful information. In general, the computation of $\pi^2(X)$ is a difficult problem, and a full discussion is well outside the scope of this modest section. We will content ourselves with sketching a few details, most of which are already contained in Taylor's paper \cite{Taylor:2012}. 

We begin by considering the following diagram:
\begin{equation}\label{pi2diag1}
\begin{gathered}
\xymatrix{S^1\ar[d]\ar[r]&S^3\ar[d]\ar[r]^-\eta&S^2\ar[d]\ar[r]^-i&\mathbb{CP}^\infty\ar@{=}[d]\ar[r]^-\rho&\mathbb{HP}^\infty\ar[d]\\
\Omega \mathbb{CP}^\infty\ar[r] &\Omega E_6\ar[r]&F\ar[r]&\mathbb{CP}^\infty\ar[r]^-{\widetilde\rho}&E_6.}
\end{gathered}
\end{equation}
The top row is a well-known fibre sequence in which $i$ classifies a generator of $H^2(S^2)$ and~$\rho$ is induced by the extension of scalars $\mathbb{C}\subseteq\mathbb{H}$. The bottom row is obtained thus. We denote by $E_6$ the 6th Postnikov section of $\mathbb{HP}^\infty$, with nonvanishing homotopy groups only in degrees~$4,5,6$. The arrow $\mathbb{HP}^\infty\rightarrow E_6$ is taken to be the natural map, and $\widetilde\rho$ is defined to be the composite of $\rho$ with this map. The space $F$ denotes the homotopy fibre of $\widetilde\rho$. Thus the bottom row of \eqref{pi2diag1} is a homotopy fibration sequence. We let $S^2\rightarrow F$ be an induced map of homotopy fibres, and complete the left-most vertical arrows by looping. 

The goal is to apply the functor $[X,-]$ to the diagram and use the exactness of the rows to study $[X,S^2]=\pi^2(X)$. We leverage the following facts. Firstly, the inclusion of the bottom cell $S^4\hookrightarrow \mathbb{H}P^4$ is 7-connected, as is the map $\mathbb{H}P^\infty\rightarrow E_6$. Consequently
$$\pi^4(X)\cong [X,\mathbb{H}P^\infty]\cong [X,E_6]
\quad\text{and}\quad
\pi^3(X)\cong [X,\Omega E_6].$$
Furthermore, $S^2\rightarrow F$ is 6-connected, so $\pi^2(X)\cong [X,F]$. 

The last fact we will use is that $E_6$ is an infinite loop space. Indeed, it approximates $S^4$ in the stable range, so we may assume that it is the $N$-fold loop space on a Postnikov section approximating the $(N+6)$-type of $S^{4+N}$ for any $N\geq0$. This is particularly important because it allows us to use the following theorem of Taylor.
\begin{thm}[Taylor~\cite{Taylor:2012}, Theorem~5.2]\label{thm_taylor thm}
Let $B$ and $C$ be homotopy associative H-spaces and $p\colon E\rightarrow B$ the homotopy fibre of a map $w\colon B\rightarrow C$. Suppose that $X$ is a space. Then, for each $\alpha\in[X,B]$ such that $w_*(\alpha)=\ast\in[X,C]$, there is a group homomorphism 
\[
\psi_\alpha\colon [X,\Omega B]\rightarrow [X,\Omega C],
\]
and a bijection from $\coker(\psi_\alpha)$ onto $(p_*)^{-1}(\alpha)\subseteq[X,E]$. \qed
\end{thm}

Let $h\colon\pi^2(S^2)\rightarrow H^2(X;\Z)$ be the map defined by $h(\alpha)=\alpha^*(s_2)$, where $s_2\in H^2(S^2;\Z)$ is a fixed generator. Also let $(-)^2\colon H^2(X;\Z)\rightarrow H^4(X;\Z)$ denote the squaring map $u\mapsto u^2$. Finally, with $\rho_2$ denoting reduction mod $2$, introduce the secondary operation
\begin{equation}\label{secondsecondaryoperation}
\Theta\colon\{u\in H^2(X;\mathbb{Z})\mid \rho_2(u^2)=0\}\rightarrow H^5(X;\mathbb{Z}/2)/(Sq^2H^3(X;\mathbb{Z}/2))
\end{equation}
which is based on the Adem relation $Sq^2(Sq^2\rho_2)=0$. The indeterminacy is unwieldy, but we will only need to use it when $Sq^2H^3(X;\mathbb{Z}/2)=0$.
\begin{prop}\label{piupper2prop}
Let $X$ be a 5-dimensional CW complex. 
\begin{enumerate}
\item If $Sq^2\colon H^3(X;\mathbb{Z}/2)\rightarrow H^5(X;\mathbb{Z}/2)$ is onto, then there is an exact sequence of sets
\[
\ast\rightarrow[X,S^3]\rightarrow[X,S^2]\xrightarrow{h} H^2(X;\Z)\xrightarrow{(-)^2}H^4(X;\Z).
\]
\item If $Sq^2\colon H^3(X;\mathbb{Z}/2)\rightarrow H^5(X;\mathbb{Z}/2)$ is trivial, then there is an exact sequence of sets 
\[
\ast\rightarrow[X,S^3]\rightarrow[X,S^2]\xrightarrow{h} H^2(X;\Z)\xrightarrow{(-)^2\oplus \theta}H^4(X;\Z)\oplus H^5(X;\mathbb{Z}/2),
\]
where
\[
\theta(u)=\begin{cases}0&u^2\neq0\\
\Theta(u)&u^2=0.\end{cases}
\]
\end{enumerate}
\end{prop}

\begin{proof}
This is obtained in both cases by applying $[X,-]$ to the top row of \eqref{pi2diag1} to get
\[
[X,S^1]\longrightarrow
[X,S^3]\overset{\eta_*}{\longrightarrow}
[X,S^2]\overset{i_*}{\longrightarrow}
[X,\C\PP^{\infty}]\overset{\rho_*}{\longrightarrow}
[X,\mathbb{HP}^{\infty}].
\]
On the left, the fact that the fibre inclusion $S^1\rightarrow S^3$ is null-homotopic implies that the image of $[X,S^1]$ in $[X,S^3]$ contains only the constant map. Moving to the right, it is clear that $h$ is the map induced by the inclusion $i\colon S^2\hookrightarrow\mathbb{CP}^\infty$ under the identification $[X,\mathbb{CP}^\infty]\cong H^2(X)$. On the far right, we have $[X,\mathbb{HP}^\infty]\cong\pi^4(X)$, and this group is described by the short exact sequence~\eqref{steenrodexactseq}.

In case $(1)$ we have $\pi^4(X)\cong H^4(X)$. The map $\rho\colon\mathbb{CP}^\infty\rightarrow\mathbb{HP}^\infty$ classifies the bundle $\gamma_\mathbb{C}\oplus\overline\gamma_\mathbb{C}$, where $\gamma_\mathbb{C}$ is the canonical line bundle over $\mathbb{CP}^\infty$ and $\overline\gamma_\mathbb{C}$ its conjugate. The composite $\mathbb{CP}^\infty\xrightarrow{\rho}\mathbb{HP}^\infty\rightarrow K(\mathbb{Z},4)$ returns the second Chern class of $\gamma_\mathbb{C}\oplus\overline\gamma_\mathbb{C}$, which is $-c_1(\gamma_\mathbb{C})^2$. Since~$c_1(\gamma_\mathbb{C})$ generates~$H^2(\mathbb{CP}^\infty)$, we identify the squaring map. Removal of the minus sign does not affect exactness.

In case $(2)$, the sequence~\eqref{steenrodexactseq} may not split. However, we do not need to understand the full structure of $\pi^4(X)$ to verify what is being claimed. Thus we replace it with $H^4(X)\oplus H^5(X;\mathbb{Z}/2)$. The first factor will play the same role as above, and for the exactness claim it will be sufficient to understand the action of~$\theta\colon H^2(X)\rightarrow H^5(X;\mathbb{Z}/2)$ on those $u\in H^2(X)$ satisfying $u^2=0$. We will show that if $u\in H^2(X)$, then $\rho_*(u)$ is null-homotopic if and only if~$u^2=0=\theta(u)$.

To this end, we replace $\mathbb{HP}^\infty$ with its 5th Postnikov section $E_5$. We have $[X,\mathbb{HP}^\infty]\cong[X,E_6]\cong[X,E_5]$, and the situation is as follows 
\begin{equation}\label{diagoninpi2proof}
\begin{gathered}
\xymatrix{&K(\mathbb{Z}/2,5)\ar[r]^-i&E_5\ar[d]\\
X\ar[r]^-u\ar@{-->}[ur]^-{u'}&\mathbb{CP}^\infty\ar[r]^-{-x^2}\ar[ur]^-{\rho'}&K(\mathbb{Z},4)\ar[r]^-{Sq^2}&K(\mathbb{Z}/2,6).}
\end{gathered}
\end{equation}
The map $\rho'$ is obtained by projecting $\widetilde\rho$ from diagram \eqref{pi2diag1} to $E_5$, and, as explained above, $-x^2=-c_1(\gamma_\mathbb{C})^2$ is what is obtained upon projecting $\rho'$ to $K(\mathbb{Z},4)$. Since $\rho'\circ u$ is the projection of $\rho_*(u)=\rho\circ u\in[X,\mathbb{HP}^\infty]$ to $[X,E_5]$, it is this composite which we need to understand. When $u^2=0$, $\rho'\circ u$ lifts to $u'\colon X\to K(\mathbb{Z}/2,5)$, and $Sq^2(H^3(X;\mathbb{Z}/2))=0$ means that $\rho'\circ u\simeq\ast$ if and only if $u'=0$. Clearly we may identify
\[
\theta(u)=u'\in H^5(X;\mathbb{Z}/2).
\]

To proceed from here it is convenient to identify $\mathbb{CP}^\infty\simeq K(\mathbb{Z},2)$ and regard the map $-x^2\colon K(\mathbb{Z},2)\rightarrow K(\mathbb{Z}/4)$ as an unstable cohomology operation. Let $F$ be its homotopy fibre and consider the diagram
\[
\xymatrix{&&K(\mathbb{Z}/2,5)\ar[d]^-i\\
&F\ar@{-->}[ur]^-{\widetilde{Sq}^2}\ar[d]^-j&E_5\ar[d]\\
X\ar[r]^-u\ar@{-->}[ur]^-{\widetilde u}&K(\mathbb{Z},2)\ar[r]^-{-x^2}\ar[ur]^-{\rho'}&K(\mathbb{Z},4)\ar[r]^-{Sq^2}&K(\mathbb{Z}/2,6).}
\]
The relation $Sq^2(-x^2)=-2x^3=0$ gives rise to a colifting $\widetilde{Sq}^2\colon F\rightarrow K(\mathbb{Z}/2,5)$, which is constructed so that $i\circ\widetilde{Sq}^2\simeq -\rho'\circ j$ \cite[p. 56]{harper}. Similarly, since $u^2=0$, we have $\widetilde u\colon X\rightarrow F$ lifting $u$. Because $\widetilde{Sq}^2\circ\widetilde u$ has order $2$, we have $i\circ\widetilde{Sq}^2\circ\widetilde u\simeq \rho'\circ u$. Hence
\[
\theta(u)=\widetilde{Sq}^2\circ\widetilde u.
\]

On the other hand, the null composition $Sq^2\circ(-x^2)\simeq\ast$ gives rise to a secondary operation
\[
\Theta'\colon\{u\in H^2(X)\mid u^2=0\}\rightarrow H^5(X;\mathbb{Z}/2)
\]
whose indeterminacy vanishes. This acts on its domain as $\Theta'(u)=\widetilde{Sq}^2\circ\widetilde u$, so
\[
\Theta'(u)=\theta(u).
\]
However, for $u\in H^2(X)$ we have $Sq^2(\rho_2(u))=\rho_2(u^2)$, meaning that $\Theta(u)$ is defined whenever $\Theta'(u)$ is. Clearly $\Theta(u)=\Theta'(u)$ in this case.
\end{proof}

\begin{remark}
For~$u\in H^2(X;\Z)$ satisfying $u^2=0$, consider the sequence
\[
X\xrightarrow{u}\mathbb{CP}^\infty\xrightarrow{x^2}K(\mathbb{Z},4)\xrightarrow{Sq^2}K(\mathbb{Z}/2,6).
\]
Each pair of consecutive arrows composes to a null map, so the functional cohomology operation $$Sq^2_u\colon H^4(\mathbb{CP}^\infty;\Z)\rightarrow H^5(X;\mathbb{Z}/2)/Sq^2(H^3(X;\mathbb{Z}/2))$$ is defined \cite[Chapter 16.1]{MosherTangora:1968}. Inspecting~\eqref{diagoninpi2proof}, we see that $\theta(u)=Sq^2_u(x^2)\in H^5(X;\mathbb{Z}/2)$.
\end{remark}

Combining Proposition~\ref{piupper2prop} and Taylor's Theorem~\ref{thm_taylor thm}, standard methods now yield the following.
\begin{prop}\label{piupper2thm}
Let $X$ be a 5-dimensional CW complex. 
Then
\begin{enumerate}
\item $\eta_*\colon\pi^3(X)\rightarrow \pi^2(X)$ is injective with image $h^{-1}(0)=\{\alpha\in\pi^2(X)\mid \alpha^*(s_2)=0\}$ where $s_2\in H^2(S^2;\Z)$ is a generator.
\item If $Sq^2\colon H^3(X;\mathbb{Z}/2)\rightarrow H^5(X;\mathbb{Z}/2)$ is onto, then there is a pairwise-disjoint decomposition \[
\pi^2(X)=\bigcup \{h^{-1}(u)\mid u\in H^2(X),\;u^2=0\}.
\]
If $u\in H^2(X;\Z)$ satisfies $u^2=0$, then $h^{-1}(u)\subseteq\pi^2(X)$ is nonempty, and $h^{-1}(u)\cong \pi^3(X)/\psi_u(H^1(X))$, where $\psi_{u}$ is the group homomorphism given in Theorem~\ref{thm_taylor thm}.
\item If $Sq^2\colon H^3(X;\mathbb{Z}/2)\rightarrow H^5(X;\mathbb{Z}/2)$ is trivial, then there is a pairwise-disjoint decomposition
\[\pi^2(X)=\bigcup\{h^{-1}(u)\mid u\in H^2(X),\;u^2=0=\Theta(u)\}.
\]
If $u\in H^2(X;\Z)$ and $u^2=0=\Theta(u)$, then $h^{-1}(u)\subseteq\pi^2(X)$ is nonempty and $h^{-1}(u)\cong \pi^3(X)/\psi_u(H^1(X))$, where $\psi_{u}$ is the group homomorphism given in Theorem~\ref{thm_taylor thm}.
\end{enumerate}
\end{prop}

\begin{remark}
For $X$ as in Proposition~\ref{piupper2thm} and $u\in H^2(X;\Z)$ lifting to $\widetilde u\colon X\rightarrow S^2$, there is an explicit description of the homomorphism $\psi_u$ given in Taylor~\cite[Section 5]{Taylor:2012}. Since~$\Omega\rho\simeq\ast$, this gives for any $\alpha\in H^1(S;\Z)$ that the map $\psi_u(\alpha)$ is the composite
\[
X\xrightarrow{\overline\Delta}X\wedge X\xrightarrow{\alpha\wedge \widetilde u}S^1\wedge S^2\xrightarrow{\phi}S^3,
\]
for some map $\phi$. The map $\phi$ must extend over $S^1\wedge\mathbb{C}P^2$ and hence must have even degree. In fact, according to \cite[Corollary 5.5, Lemma 6.5]{Taylor:2012}, $\phi$ has degree $\pm2$. 
\end{remark}
We have two applications for Proposition~\ref{piupper2thm}, which we give below. In each case it shows that the set $\pi^2(X)$ is determined by stable data.
\begin{cor}\label{piupper2thmco1}
If $X$ is an orientable, 5-dimensional Poincar\'e duality complex, $H_1(X;\Z)$ is torsion-free, and $H_2(X;\Z)$ is torsion, then $\eta\colon S^3\rightarrow S^2$ induces a bijection $\pi^2(X)\cong \pi^3(X)$. 
\end{cor}
\begin{proof}
Because $H_2(X)$ is torsion, duality gives $H^2(X)\cong H_3(X)=0$ (compare table~\eqref{table_original M hmlgy}). Thus the statement follows from part~$(1)$ of Proposition~\ref{piupper2thm}.
\end{proof}
We remark that a class of smooth, non-simply connected 5-manifolds satisfying the assumption of Corollary~\ref{piupper2thmco1} is studied in~\cite[Theorem 1.3]{KreckSu:2017}.
\begin{cor}\label{piupper2thmco2}
Let $X$ a simply connected, orientable, 5-dimensional Poincar\'e duality complex. The following statements hold.
\begin{enumerate}
\item If $Sq^2\colon H^3(X;\mathbb{Z}/2)\rightarrow H^5(X;\mathbb{Z}/2)$ is onto, then
\[
\pi^2(X)=\pi^3(X)\times H^2(X;\mathbb{Z})
\]
\item If $Sq^2\colon H^3(X;\mathbb{Z}/2)\rightarrow H^5(X;\mathbb{Z}/2)$ is trivial, then
\[
\pi^2(X)\cong \pi^3(X)\times \{x\in H^2(X;\mathbb{Z})\mid \Theta(x)=0\}
\]
where $\Theta:H^2(X;\mathbb{Z})\rightarrow H^5(X;\mathbb{Z}/2)$ is the secondary operation~\eqref{secondsecondaryoperation}.
\end{enumerate}
\end{cor}
\begin{proof}
The hypothesis of simple connectivity gives $H^1(X)=0$, meaning that each of the $\psi_u$ homomorphisms is trivial. Since $H^4(X)\cong H_1(X)=0$, the condition $u^2=0$ on $u\in H^2(X)$ is vacuous. Thus in light of Proposition~\ref{piupper2thm} the result is clear.
\end{proof}

\begin{proof}[Proof of Theorem~\ref{taylorespi2xthrm}]
The first statement is implied by Corollary~\ref{piupper2thmco2}. When $X$ is nonspin, part $(1)$ of this corollary applies and the result is immediate. Thus we reduce to the case in which $X$ is spin. Then the secondary operation $\Theta:H^2(X;\mathbb{Z})\rightarrow H^5(X;\mathbb{Z}/2)$ is defined without indeterminacy. It is known that $\Theta$ detects $\eta^2$, and, as has been discussed in Corollary~\ref{corollarysimplyconnectclass}, Barden's classification result (Theorem~\ref{Bardenclass}) implies that for this reason $\Theta$ vanishes on the cohomology of any simply connected 5-manifold.

The second statement of the theorem is a special case of Corollary~\ref{piupper2thmco1}.
\end{proof}



\begin{thebibliography}{99}

\bibitem{ACampoKotschick:1994}
N. A'Campo, D. Kotschick, {\it Contact structures, foliations, and the fundamental group},	Bull. London Math. Soc. \textbf{26}(1) (1994), 102--106.


\bibitem{arkowitz}
M. Arkowitz, {\it Introduction to homotopy theory}, Springer (2011). 


\bibitem{baues}
H. Baues, {\it Homotopy type and homology}, Oxford Mathematical Monographs, Clarendon Press, Oxford University Press (1996). 

\bibitem{BH}
H. Baues, M. Hennes, {\it The homotopy classification of $(n-1)$-connected $(n+3)$-dimensional polyhedra, $n\geq 4$}, Topology \textbf{30} (1991), 373--408. 

\bibitem{Barden:1965}
D. Barden, {\it Simply connected five-manifolds}, Ann. Math. \textbf{82} (1965), 365--385.

\bibitem{Borsuk:1936}
K. Borsuk, {\it Sur les groupes des classes de transformations 
continues}, C. R. Acad. Sci. Paris \textbf{202} (1936), 1400--1403. 


\bibitem{ChoiKajiTheriault:2017}
S. Choi, S. Kaji, S. Theriault, {\it Homotopy decomposition of a suspended real toric space}, Bol. Soc. Mat. Mex. \textbf{23} (2017), 153--161.

\bibitem{CS}
T. Cutler, T. So, {\it The homotopy type of a once-suspended 6-manifold and its applications}, Top. and App. \textbf{318} (2022).

\bibitem{GitlerStasheff:1965}
S. Gitler, J. Stasheff,	{\it The first exotic class of BF}, Topology \textbf{4} (1965), 257--266. 


\bibitem{HambletonSu:2013}
I. Hambleton, Y. Su, {\it On certain 5-manifolds with fundamental group of order 2}, Quart. J. Math. \textbf{64} (2013), 149--175.

\bibitem{HasuKishimotoSato:2016}
S. Hasui, D. Kishimoto, T. Sato, {\it $p$-local stable splitting of quasitoric manifolds}, Osaka J. Math. \textbf{53} (2016), 843--854.


\bibitem{harper}
J. Harper, {\it Secondary Cohomology Operations}, Graduate Studies in Mathematics, \textbf{49}. American Mathematical Society, Providence, RI, (2002). 

\bibitem{Huang:2021}
R. Huang, {\it Homotopy of gauge groups over non-simply connected five-dimensional manifolds}, Sci. China Math. \textbf{64} (2021), 1061--1092.

\bibitem{Huang:2022}
R. Huang, {\it Homotopy of gauge groups over high-dimensional manifolds}, Proc. Roy. Soc. Edinburgh Sect. A \textbf{152} (2022), 182--208.

\bibitem{Huang:2023}
R. Huang, {\it Suspension homotopy of 6–manifolds}, Alg. Geo. Top. \textbf{23} (2023), 439--460. 

\bibitem{HuangLi:2023}
R. Huang, P. Li, {\it Suspension homotopy of simply connected 7-manifolds}, preprint, arXiv:2208.13145.

\bibitem{JamesWhitehead:1954}
I. James, J. Whitehead,	{\it The homotopy theory of sphere bundles over spheres I}, Proc. London Math. Soc. \textbf{4} (1954), 196--218.

\bibitem{KajiTheriault:2019}
S. Kaji, S. Theriault, {\it Suspension splittings and self-maps of flag manifolds}, Acta Mathematica Sinica, \textbf{35} (2019), 445--462.

\bibitem{KirbyMelvieichner:2012}
R. Kirby, P. Melvin, P. Teichner, {\it Cohomotopy sets of 4-manifolds}, in Proc. of the Freedman Fest, Geom. and Top. Monographs \textbf{18} (2012), 161--190.

\bibitem{Konstantis:2021}
P. Konstantis, {\it Vector bundles and cohomotopies of spin 5-manifolds}, Homology, Homotopy and App., \textbf{23}(1) (2021), 143--158.

\bibitem{KreckSu:2017}
M. Kreck, Y. Su, {\it On 5-Manifolds with free fundamental group and simple boundary links in $S^5$}, Geom. and Top. \textbf{21} (2017), 2989--3008.

\bibitem{Li:2023}
P. Li, {\it Homotopy types of suspended 4-manifolds}, to appear in Algebraic and Geometric Topology, arXiv:2211.12741.

\bibitem{LiPanWu:2023}
P. Li, J. Pan, J. Wu, {\it On modular cohomotopy groups}, Isr. J. Math. \textbf{253} (2023), 887--915.


\bibitem{MadsenMilgram:1979}
I. Madsen, R. Milgram, {\it The classifying spaces for surgery and cobordism of manifolds}, Ann. of Math. Studies, No.92, Princeton Univ. Press, Princeton, (1979).

\bibitem{McCleary:2001}
J. McCleary, {\it A user’s guide to spectral sequences}, Second Edition, Cambridge Studies in Advanced Mathematics \textbf{58}, Cambridge University Press, Cambridge, (2001).

\bibitem{MosherTangora:1968}
R. Mosher, M. Tangora, {\it Cohomology operations and applications to homotopy theory}, Harper and Row, New York, (1968).

\bibitem{MS}
J. Milnor and J. Stasheff, {\it Characteristic classes}, Princeton University Press (1974).



\bibitem{Peterson:1956}
F. Peterson, {\it Generalized cohomotopy groups}, Amer. J. Math. \textbf{78} (1956), 259--281 

\bibitem{Smale:1962}
S. Smale, {\it On the structure of manifolds}, Amer. J. Math. \textbf{84} (1962), 387--399.

\bibitem{Spanier:1949}
E. Spanier, {\it Borsuk’s cohomotopy groups}, Ann. Math. \textbf{50} (1949), 203--245. 

\bibitem{Stocker:1982}
R. St{\"o}cker, {\it On the structure of 5-dimensional Poincare duality spaces}, Comm. Math. Helv. \textbf{57} (1982), 481--510.

\bibitem{ST19}
T. So, S. Theriault, {\it The suspension of a 4-manifold and its applications}, accepted by Israel Journal of Mathematics, arXiv:1909.11129, (2019).

\bibitem{Steenrod:1947}
N. Steenrod, {\it Products of cocycles and extensions of mappings}, Ann. Math. \textbf{48}(2) (1947), 290--320.

\bibitem{Taylor:2012}
L. Taylor, {\it The principal fibration sequence and the second cohomotopy set}, in Proceedings of the Freedman Fest, Geometry and Topology Monographs \textbf{18}, Geometry and Topology, Coventry, 2012, 235--251. 


\bibitem{Whitehead:1978}
G. Whitehead, {\it Elements of homotopy theory}, Springer-Verlag New York (1978).

\bibitem{Wu:2003}
J. Wu,	{\it Homotopy theory of the suspensions of the projective plane}, Mem. Amer. Math. Soc. \textbf{162} (2003).

\bibitem{Zhubr:2001}
A. Zhubr, {\it On a paper of Barden}, J. Math. Sci. \textbf{119} (2004), 35--44.

\end{thebibliography}
\end{document}